\documentclass[
preprintnumbers,superscriptaddress,reprint,floatfix]{revtex4-2}

\usepackage[utf8]{inputenc}
\usepackage[english]{babel}

\usepackage{dcolumn}
\usepackage{amsmath,amssymb,amsthm,mathrsfs, bm}
\usepackage{amsfonts}
\usepackage{graphicx}

\title{Bibliography management: \texttt{natbib} package}

\usepackage{xcolor}
\newcommand{\edit}[1]{{\color{black}{#1}}}

\begin{document}

\title{
Transitions from Monotonicity to Chaos in Gas Mixture Dynamics \\ in Pipeline Networks
}

\author{Luke S. Baker}
\affiliation{School of Mathematical and Statistical Sciences, Arizona State University, Tempe, AZ}
\affiliation{Applied Mathematics and Plasma Physics, Los Alamos National Laboratory, Los Alamos, NM}
\email[Electronic address: ]{\{lsbaker,skazi,azlotnik\}@lanl.gov}
\author{Saif R. Kazi}
\affiliation{Applied Mathematics and Plasma Physics, Los Alamos National Laboratory, Los Alamos, NM}
\affiliation{Center for Nonlinear Studies, Los Alamos National Laboratory, Los Alamos, NM}
\author{Anatoly Zlotnik}
\affiliation{Applied Mathematics and Plasma Physics, Los Alamos National Laboratory, Los Alamos, NM}



\date{\today}
 
\begin{abstract}
The blending of hydrogen generated using clean energy into natural gas pipeline networks is proposed in order to utilize existing energy systems for their planned lifetimes while reducing their reliance on fossil fuels.  We formulate a system of partial differential equations (PDEs) that govern the flow dynamics of mixtures of gases in pipeline networks under the influence of time-varying compressor and regulator control actions. The formulation is derived for general gas networks that can inject or withdraw arbitrary time-varying mixtures of gases into or from the network at arbitrarily specified nodes.  The PDE formulation is discretized in space to form a nonlinear control system that is used to prove that homogeneous mixtures are well-behaved and heterogeneous mixtures may be ill-behaved in the sense of monotone-ordering of solutions.  \edit{We use numerical simulations to compute interfaces in the parameter region of sinusoidal boundary conditions that delimit monotonic, periodic, and chaotic system responses.  The interfaces suggest that any solution in the monotonic response region is not chaotic and will eventually approach a periodic orbit.  The results are demonstrated using examples for a single pipeline and a small test network.}
\end{abstract}

\maketitle

\section{Introduction}


Although natural gas is projected to be a primary fuel source through the year 2050 \cite{nalley2022annual}, societies worldwide are investing intensively to transition from fossil fuels such as natural gas and coal to more sustainable and cleaner resources.  Hydrogen is an energy carrier that can be cleanly produced \cite{salvi2015sustainable} and can address climate change because it does not result in carbon dioxide emissions or other harmful emissions when it is burned.   Several qualities of hydrogen make it an attractive fuel option for a variety of applications that include transportation and high temperature manufacturing.  Hydrogen can also be used to power turbines, which can potentially be used for aviation and electric power production.  Hydrogen can be produced directly from fossil fuels, biomass, or direct electrolysis, by splitting water into its constituent components of hydrogen and oxygen.   After hydrogen is produced, it can be transported to end users economically by dedicated pipeline systems.


Recent studies have proposed that natural gas pipelines can safely transport mixtures of up to 20\% hydrogen or more by volume \cite{gotz2016renewable, ozturk2021comprehensive}.  Hydrogen could be transported through the existing infrastructure and then separated, or the mixture could be burned directly as an end-use fuel.  Because the physical and chemical properties of hydrogen and natural gas (primarily methane) differ significantly, the mass and energy transport dynamics of inhomogeneous mixtures of these constituent gases are considerably more complex than for a homogeneous gas \cite{van2004math}.  The mathematical modeling of such mixtures is also considerably more challenging than what has traditionally been done for gas pipelines \cite{melaina2013blending}.    
The introduction of substantial proportions of much lighter hydrogen into natural gas pipelines requires much closer spacing of gas compressors, and this relationship has been characterized in an empirical study \cite{WITKOWSKI20172508}.  Additionally, the pressure and flow dynamics in gas networks have been proven to satisfy certain physically intuitive and conceptually valuable monotonicity properties \cite{misra2020monotonicity}, which must be re-examined in the presence of inhomogeneous gas mixing.


The physical complexities of blending hydrogen in natural gas pipelines present several mathematical challenges.  Additional state variables are needed to account for changes in mass fraction, which affect total density, energy content, and flow dynamics.  Modeling the flow of a homogeneous gas on a network requires partial differential equations (PDEs) for mass and momentum conservation on each pipe, and a linear mass flow balance equation at each network junction.  Injection of a second gas into the network requires the addition of another PDE on each pipe and a bilinear nodal balance equation at each junction to account for conservation of composition.  This more than doubles the state space of the continuous model.  Moreover, the faster wave speed corresponding to the lower density of hydrogen worsens the numerical ill-conditioning of the dynamic model.  Such issues have been highlighted by the numerical simulations of hydrogen and natural gas flows in pipelines \cite{uilhoorn2009dynamic, chaczykowski2018gas, guandalini2017dynamic, elaoud2017numerical, hafsi2019computational, fan2021transient, agaie2017reduced, WITKOWSKI20172508, SUBANI2017244}. 

One recent study has demonstrated conditions under which pipeline pressures may exceed allowable upper limits, and that the likelihood of this occurrence increases proportionally with increasing hydrogen concentration \cite{hafsi2019computational}.  Another study examined the effects of hydrogen blending on the detection and estimation of leaks \cite{SUBANI2017244}, and demonstrated that the amount of leak discharge increases as the concentration of hydrogen increases.  A moving grid method and an implicit backward difference method for tracking gas concentration were both shown to perform well for numerical simulations but the implicit difference method may lose some finer detail due to numerical diffusion \cite{chaczykowski2018gas}.   The method of characteristics was also applied for the numerical simulation of transient flows on cyclic networks with homogeneous flow mixtures \cite{elaoud2017numerical}.  Modeling networks of pipelines with composition tracking was the focus of another recent study \cite{WITKOWSKI20172508}, although this model does not include control actions of compressor units.  In general, these models demonstrate a simulation capability or sensitivity study for a specific network.  Addressing challenging design, operational, and economic issues in the pipeline transport of gas mixtures will require minimal and generalizable mathematical models that adequately describe the relevant physics, in addition to complex simulations with comprehensive characterizations of pipeline flows.


The scope of the present study is threefold.  First, we extend general control system models for gas pipeline networks \cite{zlotnik2015optimal} to account for heterogeneous mixtures of hydrogen and natural gas.  \edit{Similar lumped parameter modeling has been well-studied for pipeline simulation since at least a decade ago \cite{grundel2013computing,grundel2014model}.} The state variables are flows, partial densities, and pressures throughout the network, and the control variables are the actions of compressor and regulator units.  Control actions may be designed to minimize fuel consumption \cite{wong1968optimization, percell1987steady, rachford2000optimizing} or maximize economic value \cite{zlotnik2019scheduling}.  The PDE control system of the mixture is discretized in space using a finite volume method \cite{himpe2021model} and written in matrix form as a finite-dimensional control system of nonlinear ordinary differential equations (ODEs). 
Second, we prove that solutions to initial boundary value problems (IBVPs) of a gas mixture have certain monotone ordering properties if the concentration is homogeneous but, in general, do not have these properties if the concentration is heterogeneous.  The homogeneous monotonicity result generalizes the pure natural gas monotonicity result for obtaining control formulations that are robust to uncertainty in pressure and withdrawal profiles \cite{zlotnik2016monotonicity, misra2020monotonicity}.  \edit{Third, we demonstrate that the solution of an IBVP may be irregular, in the sense of generating a continuous distribution of harmonic modes, and may also be chaotic, in the sense of being sensitive to initial conditions.  Numerical simulations are used to characterize flow solution behavior in a phase space of periodic forcing functions and to identify boundaries between the regions of monotonic, periodic, and chaotic solution behavior.  Transitions through such fluid mixing phase regions were observed in oceanic wind bursts \cite{tziperman1994nino, eisenman2005westerly} and in flame combustion of hydrogen and air mixtures \cite{pizza2008dynamics, alipoor2016combustion}.  Inspection of the response interfaces that we compute suggest that any solution in the monotonic response region is not chaotic and will eventually approach a periodic orbit.} Simulation-based analyses such as those presented here could be used to evaluate appropriate limitations on blending of hydrogen into existing natural gas pipeline networks.

The rest of this paper is organized as follows.  The PDEs that govern heterogeneous mixtures of hydrogen and natural gas are presented in Section \ref{sec:pde_flow}.  In Section \ref{sec:discretize_flow}, the PDE system is discretized in space to obtain a system of ODEs.  Section \ref{sec:equivalent_sys} presents a derivation of equivalent ODE systems in terms of other state variables of interest.   Section \ref{sec:monotonicity} contains a proof that each of the equivalent systems have monotonic solutions if the concentration is homogeneous, as well as a proof that the solutions are, in general, non-monotonic if the concentration is heterogeneous.  In Section \ref{sec:network}, we illustrate non-monotonic system responses using numerical simulations of flows through a small test network that contains a loop, and which was examined in a previous study \cite{GYRYA201934}.  Moreover, that section illustrates that certain types of equivalent systems may have more desirable monotone system behavior than others in certain response regimes.  \edit{Sections \ref{sec:MI}, \ref{sec:PI}, and \ref{sec:chaotic} describe techniques to compute interfaces, in the region of boundary condition parameters for flow in a single pipeline, between regions that do and do not exhibit monotonic, periodic, and chaotic properties, respectively.  We provide concluding remarks and an outlook for future work in Section \ref{sec:conclusion}.}

 \vspace{-2ex}
\section{Gas Network Modeling}  \label{sec:pde_flow}

A gas transport network is modeled as a connected and directed graph $(\mathcal E, \mathcal V)$ consisting of edges $\mathcal E =\{1,\dots,E\}$ and nodes $\mathcal V=\{1,\dots,V\}$, where $E$ and $V$ denote the numbers of edges and nodes, respectively.  It is assumed that the elements of these sets are ordered according to their integer labels.  The edges represent pipelines and the nodes represent junctions or stations where gas can be injected into or withdrawn from the network.  The symbol $k$ is reserved for indexing edges in $\mathcal E$ and the symbols $i$ and $j$ are reserved for indexing nodes in $\mathcal V$. 
The graph is directed by assigning a positive flow direction along each edge.  It is assumed that gas physically flows in only the direction of positive flow, so that the mass flow and velocity values of the gas are positive quantities everywhere in the network.  The notation $k:i\mapsto j$ means that edge $k\in \mathcal E$ is directed from node $i\in \mathcal V$ to node $j\in \mathcal V$.  For each node $j\in \mathcal V$, we define (potentially empty) incoming and outgoing sets of pipelines by $_{\mapsto} j=\{k\in \mathcal E| k:i\mapsto j\}$ and $j_{\mapsto}=\{k\in \mathcal E| k:j\mapsto i \}$, respectively.   \edit{All nomenclature is listed in Appendix \ref{sec:nomenclature}.}

 \vspace{-1ex}
\subsection{Modeling Physical Flow in a Pipe} \label{subsec:pipeflow}

Compressible flow of a homogeneous ideal gas through a pipe is described using the one-dimensional isothermal Euler equations \cite{osiadacz1984simulation},
\begin{subequations} \label{eq:gaspde0}
\begin{align}
    \partial_t \rho + \partial_x (\rho u) & = 0, \label{eq:gaspde0a} \\
    \partial_t (\rho u) + \partial_x (p + \rho u^2) & = - \frac{\lambda}{2D}\rho u |u| - \rho g \frac{\partial h}{\partial x}, \label{eq:gaspde0b} \\
    p = \rho ZR\mathbf{T} & = \sigma^2 \rho, \label{eq:gaspde0c}
\end{align}
\end{subequations}
where the variables $u(t,x)$, $p(t,x)$, and $\rho(t,x)$ represent velocity, pressure, and density of the gas, respectively.  Here, $t\in [0,T]$ and $x\in[0,\ell]$, where $T$ denotes the time horizon and $\ell$ denotes the length of the pipe.  The symbols $\partial_t$ and $\partial_x$ denote the differential operators with respect to time $t$ and axial location $x$, respectively.  The above system describes mass conservation \eqref{eq:gaspde0a}, momentum conservation \eqref{eq:gaspde0b}, and the gas equation of state \eqref{eq:gaspde0c}.   The variable $h$ represents the elevation of the pipe.  The dominant term in the momentum equation \eqref{eq:gaspde0b} is the phenomenological Darcy-Weisbach term that models momentum loss caused by turbulent friction, and is scaled by a dimensionless parameter $\lambda$ called the friction factor.  The remaining parameters are the internal pipe diameter $D$, the wave (sound) speed  $\sigma=\sqrt{ZR\mathbf{T}}$ in the gas, and the gravitational acceleration $g$, where $Z$, $R$, and $\mathbf{T}$ are the gas compressibility factor, specific gas constant, and absolute temperature, respectively.  Here, we assume that gas pressure $p$ and gas density $\rho$ satisfy the ideal gas equation of state \eqref{eq:gaspde0c} with wave speed $\sigma$.  While non-ideal modeling is necessary in practice to correctly quantify flows at pressures used in large gas transport pipelines, ideal gas modeling still qualitatively captures the flow phenomenology, so we use it for simplicity of exposition.  Extension to non-ideal gas modeling can be made by applying appropriate nonlinear transforms \cite{GYRYA201934}. 

It is standard to use the per area mass flux $\varphi=\rho u$, and assume that gas flow is an isothermal process, that flow is turbulent and has high Reynolds number, and that the flow is adiabatic, i.e. there is no heat exchange with the ground \cite{herty2010new}.  For slowly varying boundary conditions, the kinetic energy term $\partial_x(\rho u^2)$ and the inertia term $\partial_t (\rho u)$ in equation \eqref{eq:gaspde0b} may be omitted \cite{osiadacz1984simulation}.  With these assumptions, and given no elevation changes, Eq. \eqref{eq:gaspde0} can be reduced to
\begin{subequations}
\label{eq:pde_1}
\begin{align}
 \partial_t \rho + \partial_x \varphi &= 0, \label{eq:pde_1a}\\
 \partial_x (\sigma^2\rho) & = -\frac{\lambda}{2D} \frac{\varphi |\varphi|}{\rho}. & \label{eq:pde_1b}
\end{align}
where $\rho$ and $\varphi$ denote density and mass flux (in per-area units).  The above set of equations have been used in several previous studies \cite{sundar2019tcst, misra2020monotonicity}, and we refer the reader there for further justifications.  
Here, we extend these equations to the case of a mixture of two constituent gases, whose partial pressures, partial densities, partial fluxes, and mass fractions are denoted by $p^{(m)}$, $\rho^{(m)}$, $\varphi^{(m)}$, and $\eta^{(m)}$, respectively, where $m=1$ and $m=2$ are used to identify the two distinct gases.  The fraction of mass of each gas is related to the partial density variables by $\eta^{(m)}= \rho^{(m)}/(\rho^{(1)}+\rho^{(2)})$.  The propagation of either mass fraction quantity $\eta^{(m)}$ can be modeled by the convection-diffusion equation with diffusion terms omitted \cite{chaczykowski2018gas}, i.e.,
\begin{align}
& \partial_t \eta^{(m)} + \frac{\varphi}{\rho}\partial_x \eta^{(m)} = 0. \label{eq:pde_convect}
\end{align}
\end{subequations}
It follows from the relation $\eta^{(1)}=(1-\eta^{(2)})$ that only one mass fraction variable should be modeled if there are exactly two distinct gases in the mixture.  It can be shown that $\eta^{(1)}$ solves Eq. \eqref{eq:pde_convect} if and only if $\eta^{(2)}$ does.  Partial densities $\rho^{(m)}$ and partial fluxes $\varphi^{(m)}$ are proportional to their total counterparts $\rho=(\rho^{(1)}+\rho^{(2)})$ and $\varphi=(\varphi^{(1)}+\varphi^{(2)})$ with the mass fraction being the coefficient of proportionality.  In particular, $\rho^{(m)}=\eta^{(m)}\rho$ and $\varphi^{(m)}=\eta^{(m)}\varphi$.  However, proportionality according to mass fraction does not generally hold true for pressure.  The ideal equation of state for each constituent gas is defined by $p^{(m)}=\sigma_m^2\rho^{(m)}$, where $\sigma_1$ and $\sigma_2$ are the wave speeds of the two gases.  Using Dalton's law \cite{silberberg2006chemistry}, the total pressure of the mixture is defined to be the summation of partial pressures given by 
\begin{equation*}
    p=p^{(1)}+p^{(2)}=\sigma_1^2\rho^{(1)}+\sigma_2^2\rho^{(2)}=\sigma^2\rho,
\end{equation*}
where $\sigma^2=(\sigma_1^2\eta^{(1)}+\sigma_2^2\eta^{(2)})$.  It follows that the local wave speed $\sigma$ of the mixture depends on the local mass fraction of the gases.  Moreover, the total pressure and partial pressure variables are related through the {\it volumetric} concentration defined by $\nu^{(m)}=\sigma_m^2\rho^{(m)}/(\sigma_1^2\rho^{(1)}+\sigma_2^2\rho^{(2)})$.  
In particular, $p^{(m)}=\nu^{(m)}p$.  From here onward, we use the terms mass fraction and concentration interchangeably, and specifically refer to \emph{volumetric} concentration where that quantity is examined. 

\edit{In our application of energy transport, superscripts ``(1)'' and ``(2)'' will henceforth identify correspondence of variables to natural gas and hydrogen, respectively.  The wave speeds of natural gas and hydrogen are defined, respectively, by $\sigma_1$ and $\sigma_2$.  We note that natural gas is itself a mixture composed of several distinct gases, and its composition has historically been modeled as constant and homogeneous in the majority of academic pipeline simulation studies.  This simplification is accepted because the composition of processed pipeline quality gas is at least 90\% methane and typically 95\% methane by molar mass, and with ethane comprising the majority of the other constituents \cite{kunz2012gerg}.  The calorific content of methane from different sources may vary somewhat \cite{hante2019complementarity}, and numerical methods for pipeline simulation that account for variation in gas composition are well developed \cite{chaczykowski2018gas}.  In contrast, the calorific values and molecular masses of hydrogen and natural gas both differ by an order of magnitude.  
In this study, we suppose that natural gas has a nominal homogeneous and constant composition, and focus on the phenomenology of dynamic pipeline response to injection of much less dense hydrogen gas with variation in time and by location. 
In general, although our focus here is on blending two gases in a network of pipelines, the models and theory here can be extended to mixtures of more than two gases with additional equations and variables corresponding to superscripts ``($m$)" for $m=1,\dots,M$, where $M$ would denote the total number of gas components of the mixture.}


\subsection{Gas Mixture Dynamics on a Network} \label{subsec:networkflow}

With the above assumptions, the flow dynamics through the horizontal pipeline of index $k\in \mathcal E$ is modeled with the friction-dominated PDEs  
\begin{eqnarray}
\partial _t \rho_k^{(m)} +\partial_x \left(\frac{\rho_k^{(m)}}{\rho_k^{(1)}+\rho_k^{(2)}} \varphi_k \right) &=&0, \label{eq:pde1} \\
\partial_x\left(\sigma^2_1\rho_k^{(1)}+\sigma^2_2\rho_k^{(2)} \right) &=& -\frac{\lambda_k}{2D_k}\frac{ \varphi_k|  \varphi_k|}{\rho_k^{(1)}+\rho_k^{(2)}}, \label{eq:pde2} \quad
\end{eqnarray}
where Eq. \eqref{eq:pde1} is defined for both $m=1$ and $m=2$.  We leave it to the reader to verify that Eqs. \eqref{eq:pde1}-\eqref{eq:pde2} defined in terms of partial densities $\rho^{(1)}$ and $\rho^{(2)}$ and total flow $\varphi$ are equivalent to Eq. \eqref{eq:pde_1}  defined in terms of total density $\rho$, total flow $\varphi$, and one concentration variable $\eta^{(2)}$.   The wave speeds, $\sigma_1$ and $\sigma_2$, and the diameter $D_k$, length $\ell_k$, and friction factor $\lambda_k$ associated with each pipeline $k\in \mathcal E$ are the parameters of the system.  

Compressor and regulator stations are critical components that actuate the flow of gas through the network and reduce pressure in the direction of flow, respectively.  For convenience, we assume that a compressor is located at the inlet and a regulator is located at the outlet of each pipeline, where inlet and outlet are defined with respect to the oriented positive flow direction.  For each pipeline $k\in \mathcal E$, compression and regulation are modeled with multiplicative control variables $\underline\mu_k(t)\ge1$ and $\overline \mu_k(t)\ge1$, respectively.  That is, discharge pressure and density are $\underline\mu_k(t)$ times larger than the suction pressure and density.  

\begin{figure}[!t]
\centering
\includegraphics[width=.75\linewidth]{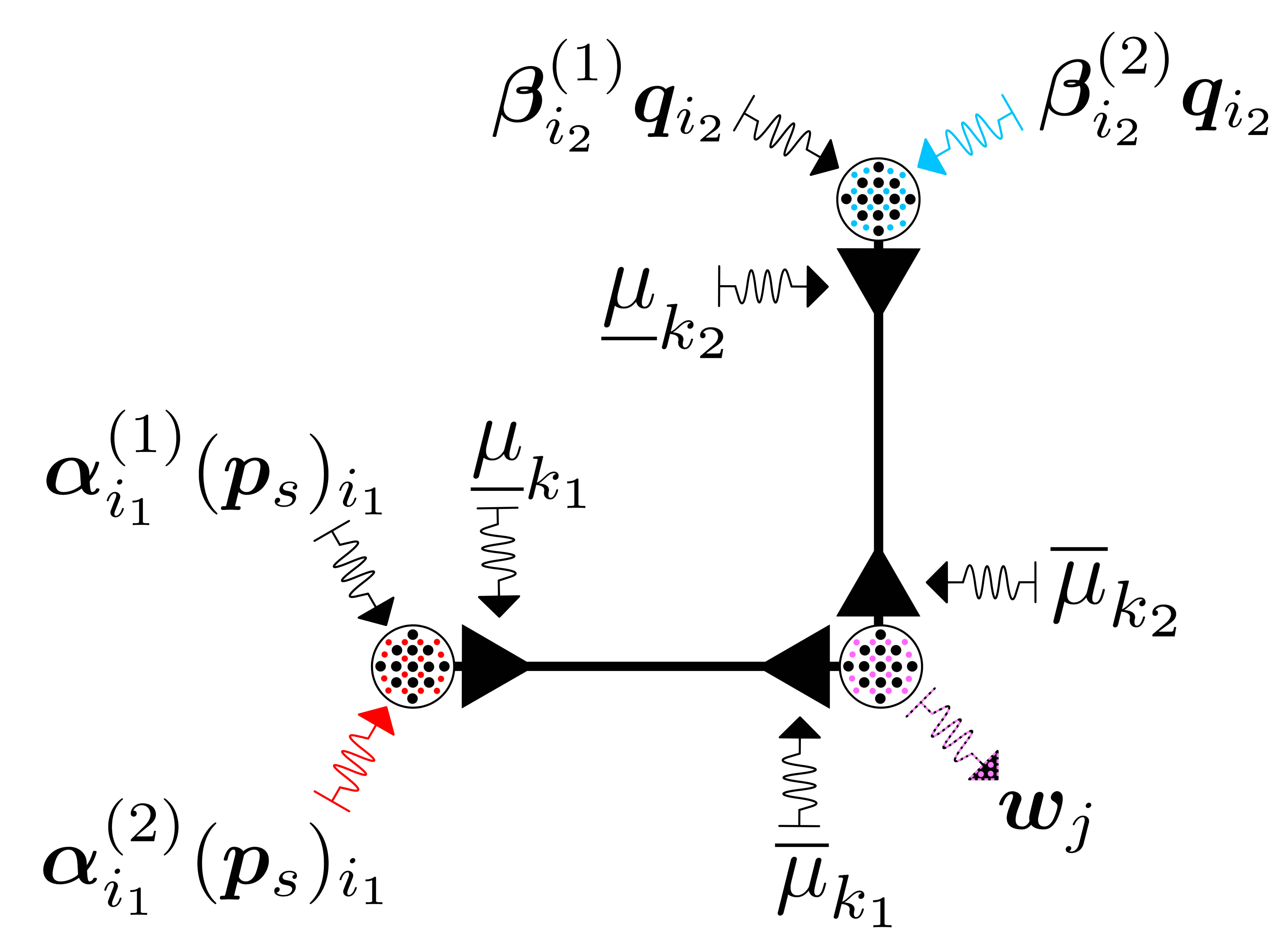}
\caption{Configuration of the boundary conditions.  Here,  $j\in \mathcal V_w$, $k_1:i_1\mapsto j$ with $i_1\in \mathcal V_s$, and $k_2:i_2\mapsto j$ with $i_2\in \mathcal V_q$.}
\label{bc_conf}
\end{figure}

The boundary conditions for a mixture of gases allow for more degrees of freedom than those for a single gas, and are formulated here to enable definition of a range of potential scenarios.  All of the flow quantities defined in this paragraph are, in general, time-varying, but we suppress time-dependence for readability.  The network nodes are partitioned into slack nodes $ \mathcal V_s\subset \mathcal V$ and non-slack nodes $ \mathcal V_d\subset \mathcal V$.  \edit{Slack (pressure) nodes are typically used to represent large sources of gas, at which pressure and concentration are specified and inflow to the network is a dependent variable that is determined by solving the network flow equations defined below.  At non-slack (flow) nodes, the withdrawal flow from the network is specified, and the pressure and mass fraction is determined by solving the initial boundary value problem.  Alternatively, injection can be specified at a non-slack (flow) node if concentration is also provided.}  Slack nodes are assumed to be ordered in $\mathcal V$ before non-slack nodes, so that $i<j$ for all $i\in \mathcal V_s$ and $j\in \mathcal V_d$.  A mixture of gas is injected into the network at each slack node $i\in \mathcal V_s$.  The boundary conditions at the slack nodes $i\in \mathcal V_s$ are defined by specifying individual densities $ \bm s_i^{(1)}$ and $\bm s_i^{(2)}$.  Alternatively, pressure $(\bm p_s)_i$ and concentration $ \bm \alpha_i^{(m)}$ may be specified at slack nodes $i\in \mathcal V_s$. The relations $(\bm p_s)_i=(\sigma_1^2 \bm s_i^{(1)}+\sigma_2^2 \bm s_i^{(2)})$ and $ \bm \alpha_i^{(m)}=\bm s_i^{(m)}/(\bm s_i^{(1)}+\bm s_i^{(2)})$ can then be used to determine the corresponding partial densities that will achieve the specified pressures and concentrations.  Non-slack nodes are partitioned into injection nodes $\mathcal V_q\subset \mathcal V_d$ and withdrawal nodes $\mathcal V_w\subset \mathcal V_d$.  
We order the non-slack nodes $\mathcal V_d$ with injection nodes enumerated before withdrawal nodes, so that $i<j$ for all $i\in \mathcal V_q$ and $j\in \mathcal V_w$. A mixture is withdrawn from the network at each withdrawal node $j\in \mathcal V_w$ with boundary conditions specified by mass outflow $\bm w_j\ge 0$.  At each injection node $j\in \mathcal V_q$, a mixture is injected into the network with boundary conditions specified by both the mass inflow  $\bm q_j$, with $\bm q_j \ge 0$, and the concentration $\bm \beta_j^{(m)}$.   Although a mass inflow is specified at each injection node $j\in \mathcal V_q$ with concentration $\bm \beta_j^{(m)}$, this does not, in general, imply that the concentration flowing from node $j$ into outgoing edges is equal to $\bm \beta_j^{(m)}$, because the nodal concentration is a mixture of flows entering node $j$ either by injection or from incoming pipelines.  Boundary condition designations are illustrated for a small example network in Fig. \ref{bc_conf}.

For each non-slack node $j\in \mathcal V_d$, the nodal partial density variables are denoted by $\bm \rho_j^{(m)}$ and the dependent concentration variables are denoted by $\bm \eta_j^{(m)}=\bm \rho_j^{(m)}/(\bm \rho_j^{(1)}+\bm \rho_j^{(2)})$ for $m=1$ and $m=2$. All of the nodal quantities in this study are identified with bold symbols.  Inlet and outlet edge variables are defined by attaching underlines below and overlines above the associated edge variables, respectively.  For example, $\underline \varphi_k(t)= \varphi_k(t,0)$ and $\overline \varphi_k(t)= \varphi_k(t,\ell_k)$.  Let us denote the cross-sectional area of edge $k\in \mathcal E$ by $\chi_k=\pi D_k^2/4$.  The boundary conditions for the flow of the mixture are defined for $m=1$ and $m=2$ by
\begin{eqnarray}
\underline{\rho}_{k}^{(m)}&=&\underline\mu_{k} \bm s_{i}^{(m)},\qquad 
\overline \rho_{k}^{(m)}=\overline \mu_{k} \bm \rho_{j}^{(m)}, \label{eq:pde_bc1} \\
\underline{\rho}_{k}^{(m)}&=&\underline\mu_{k} \bm \rho_{i}^{(m)},\qquad 
\overline \rho_{k}^{(m)}=\overline \mu_{k} \bm \rho_{j}^{(m)}, \label{eq:pde_bc2} \\
\bm \gamma_j^{(m)}\bm d_j&=& \sum_{k\in _{\mapsto}j}  \chi_k \overline \eta_k^{(m)} \overline\varphi_{k}
-\sum_{k\in j_{\mapsto}} \chi_k\underline \eta_k^{(m)} \underline\varphi_{k},  \label{eq:pde_bc}
\end{eqnarray}
where Eq. \eqref{eq:pde_bc1} is defined for $k:i\mapsto j$ with $i\in \mathcal V_s$, Eq. \eqref{eq:pde_bc2} is defined for  $k:i\mapsto j$ with $i, j\in \mathcal V_d$, and Eq. \eqref{eq:pde_bc} is defined for $j\in \mathcal V_d$ with the condition that $\bm \gamma_j^{(m)}\bm d_j=\bm \eta_j^{(m)}\bm w_j$ if $j\in \mathcal V_w$ and $\bm \gamma_j^{(m)}\bm d_j=-\bm \beta_j^{(m)}\bm q_j$ if $j\in \mathcal V_q$. \edit{The conditions in  Eqs. \eqref{eq:pde_bc1}-\eqref{eq:pde_bc2} are multiplicative relations between the partial densities at the boundaries of pipelines and the auxiliary ``nodal'' variables for partial densities that are internal to the nodes at each end (see Fig. \ref{bc_conf}).  The multiplicative factors are used as control variables that represent gas compression and regulation.  The condition for the conservation of mass through each non-slack node in Eq. \eqref{eq:pde_bc} depends on whether the node is an injection node or a withdrawal node. For injection, we specify both the concentration of the mixture and the mass flow being injected into the node, but for withdrawal, we specify only the mass outflow. The concentration of the withdrawn gas is determined by solving the network flow equations.}

The initial conditions of partial density are assumed to be a steady-state solution given for all $k\in \mathcal E$ and $x\in [0,\ell_k]$ by
\begin{eqnarray}
 \rho_k^{(m)}(0,x)=\varrho_k^{(m)}(x). \label{eq:pde_ic}
\end{eqnarray}
The steady-state configuration is defined to be the solution of the system in Eqs. \eqref{eq:pde1}-\eqref{eq:pde_bc} when the boundary condition profiles are time-invariant (i.e. equal to the initial values of the time-varying boundary profiles).  More details on the initial condition for the discretized system are provided in the following section.  Mass flux is not initially specified because it is uniquely determined from the density.  We assume standard conditions for well-posedness \cite{gugat2012well}, and specifically that the boundary conditions are smooth, slowly-varying, bounded in their respective domains, and compatible with the initial conditions to ensure the existence of a smooth, slowly-varying, bounded solution.   The flow of the mixture of gases in the network is defined by the initial-boundary value system of PDEs defined by equations \eqref{eq:pde1}-\eqref{eq:pde_ic}.

\section{Spatial Discretization} \label{sec:discretize_flow}

To analyze the system of PDEs \eqref{eq:pde1}-\eqref{eq:pde_ic} on the graph $(\mathcal E, \mathcal V)$, we have developed a process of discretization, which includes a \emph{refinement} of the graph, \emph{approximation} of the PDE system by an ODE system using a finite volume approach, and a \emph{reformulation} in terms of variable vectors and parameter matrices.  The vectors include variables that represent the states and boundary condition profiles, and the matrices incorporate network model parameters, the incidence structure of the graph, and the control values.

{\bf Graph Refinement:} A refinement $(\hat{\mathcal E},\hat{\mathcal V})$ of the graph $(\mathcal E, \mathcal V)$ is created by adding auxiliary nodes to $\mathcal V$ in order to subdivide the edges of $\mathcal E$ so that $\ell_k\le \ell$ for all $k\in \hat{\mathcal{E}}$, where $\ell$ is sufficiently small \cite{grundel2013computing}.  Henceforth, we assume that $\ell\le 1$ (km), and will use that threshold for computational studies as well.  The refined graph inherits the prescribed orientation of the parent graph.  Assuming sufficiently fine network refinement, the relative difference of the density variables of adjacent nodes in the solution to the IBVP \eqref{eq:pde1}-\eqref{eq:pde_ic} can be made arbitrarily small in magnitude because of continuity of the solution to the system given well-posed conditions \cite{gugat2012well}.  We assume for all $k\in \hat{\mathcal E}$ that
\begin{equation}
     \frac{ \left| \overline \rho_k^{(m)}-\underline \rho_k^{(m)} \right|}{ \underline \rho_k^{(m)}} < \epsilon, \qquad 
     \frac{ \left| \overline \rho_k^{(m)}-\underline  \rho_k^{(m)} \right|}{\overline  \rho_k^{(m)}} < \epsilon, \label{eq:eps}
\end{equation}
where $0\le \epsilon \ll 1$.  The proofs that follow only require $\epsilon \le 1$.  We assume that the graph has been sufficiently refined to satisfy Eq. \eqref{eq:eps} and that the hats may be omitted moving forward.

{\bf Finite Volume Approximation:} The system of ODEs is obtained by integrating the dynamics in Eqs. \eqref{eq:pde1}-\eqref{eq:pde2} along the length of each refined pipeline segment so that
\begin{eqnarray*}
\int_0^{\ell}\partial _t \rho^{(m)}dx = -\int_0^{\ell}\partial_x \left(\frac{\rho^{(m)}}{\rho^{(1)}+\rho^{(2)}} \varphi \right)dx,  \\
\int_0^{\ell}\partial_x\left(\sigma^2_1\rho^{(1)}+\sigma^2_2\rho^{(2)} \right)dx
= -\frac{\lambda}{2D}\int_0^{\ell}\frac{ \varphi|  \varphi|}{\rho^{(1)}+\rho^{(2)}}dx,
\end{eqnarray*}
where edge subscripts have been removed for readability. The above integrals of space derivatives are evaluated using the fundamental theorem of calculus.  The remaining integrals are evaluated by approximating pipeline density by outlet density and pipeline flux by inlet flux.  These approximations are independent of $x$ and may be factored out of the integrals.  The above equations become
\begin{eqnarray}
\ell\dot {\overline\rho}^{(m)} &=&\underline \eta^{(m)} \underline\varphi-\overline \eta^{(m)} \overline\varphi, \label{eq:pde1_integrate}
\\
\sum_{n=1}^{2}\sigma_n^2\left( \overline\rho^{(n)}-\underline \rho^{(n)} \right) &=& -\frac{\lambda \ell}{2D}\frac{\underline \varphi \left| \underline \varphi \right|}{\overline\rho^{(1)}+\overline \rho^{(2)}},  \label{eq:pde2_integrate}
\end{eqnarray}
where a dot above a variable represents the time-derivative of the variable.  

{\bf Matrix Form:} We now write the discretized system in matrix-vector form.  Define $E\times E$ diagonal matrices $L$ and $X$ with diagonal entries $L_{kk}=\ell_{k}$ and $X_{kk}=\chi_k$.  Define the time-varying (transposed) incidence matrix $M$ of size $E \times V$ componentwise by
\begin{align}
M_{kj}&=  
\begin{cases}
\overline \mu_{k}(t), & \text{ edge $k\in _{\mapsto}j$ enters node $j$},
 \\
-\underline\mu_{k}(t), & \text{ edge $k\in j_{\mapsto}$ leaves node $j$,}
\\
0, & \text{else.}
\end{cases}
\end{align}
Define the $E\times V_s$ submatrix $M_s$ of $M$ by the removal of columns $i\in \mathcal V_d$, the $E\times (V-V_s)$ submatrix $M_d$ of $M$ by the removal of columns $i\in \mathcal V_s$, and the positive and negative parts of $M_d$ by $\overline M_d$ and $\underline M_d$ so that $M_d=(\overline M_d +\underline M_d)/2$ and $|M_d|=(\overline M_d -\underline M_d)/2$, where $V_s$ denotes the number of slack nodes and $|A|$ denotes the componentwise absolute value of a matrix $A$. Define the signed matrices $Q_d=$sign$(M_d)$, $\overline Q_d=$sign$(\overline M_d)$, $\underline Q_d=$sign$(\underline M_d)$, and similarly for $M_s$.  These signed matrices are well-defined by the lower-bound constraints on compression and regulation. \edit{The incidence matrix provides a compact representation of the network structure that can be used to specify the dynamical equations for all of the refined edges. Moreover, the signed incidence matrix is used to produce the first order finite difference approximation of the spatial derivative of each component of the vector on which it acts. The absolute values of the positive and negative parts of the weighted incidence matrix apply the multiplicative boost in compression and reduction in regulation, respectively, to the vectors on which they act.} 

We define the $V_d\times V_d$ identity matrix $I$, the $V_d\times V_q$ submatrix $I_q$ of $I$ by the removal of columns $j\in \mathcal V_w$, and the $V_d\times V_d$ matrix $I_w$ by replacing columns $j\in \mathcal V_q$ of $I$ with the zero vector. Here, $V_d$ and $V_q$ denote the numbers of non-slack nodes and non-slack injection nodes, respectively. We also define inlet and outlet edge mass flux vectors by $\underline \varphi =(\underline \varphi_1,\dots,\underline \varphi_E)^T$ and $\overline \varphi =(\overline \varphi_1,\dots,\overline \varphi_E)^T$, and similarly for inlet and outlet edge concentrations $\underline \eta$ and $\overline \eta$.  Then, let us define the vectors $\bm \rho^{(m)}=(\bm \rho_{V_s+1}^{(m)},\dots, \bm \rho_{V_d}^{(m)})^T$, $\bm \alpha^{(m)}=(\bm \alpha_1^{(m)},\dots, \bm \alpha_{V_s}^{(m)})^T$, and $\bm \beta^{(m)}=(\bm \beta^{(m)}_{V_s+1},\dots,\bm \beta^{(m)}_{V_q})^T$, where the subscripts of the vector components are indexed according to the node labels in $\mathcal V$.  Similarly, define the vectors $\bm \eta^{(m)}=(\bm \eta^{(m)}_{V_s+1},\dots,\bm \eta^{(m)}_{V_d})^T$ and $\bm d=(\bm d_{V_s+1},\dots,\bm d_{V_d})^T$.  Recall that the components of $\bm d$ are positive for those corresponding to non-slack withdrawal nodes and negative for non-slack injection nodes.  Define the function $f:\mathbb{R}^E\times \mathbb{R}^E\rightarrow \mathbb{R}^E$ component-wise for $k\in \mathcal E$ by
\begin{align}
f_{k}(y,z)&=-\text{sign}(z_k)\Lambda_k \left| y_kz_k  \right|^{1/2}, \label{eq:flux_expression}
\end{align}
where $\Lambda_k =\sqrt{2D_{k}/(\lambda_{k}\ell_k)}$.  This function is used to express $\underline \varphi$ in Eq. \eqref{eq:pde2_integrate} in terms of density and its spatial derivative so that we may eliminate flux from the dynamic equations.  Using the function in Eq. \eqref{eq:flux_expression}, the discretized flow in Eqs. \eqref{eq:pde1_integrate}-\eqref{eq:pde2_integrate} together with the boundary conditions in Eqs. \eqref{eq:pde_bc1}-\eqref{eq:pde_bc} may be expressed in matrix-vector form as
\begin{eqnarray}
\!\!\!\! L\overline M_d \dot {\bm \rho}^{(m)} &=& \underline \eta^{(m)} \odot F -\overline \eta^{(m)}  \odot \overline \varphi ,  \label{eq:reduced_partial}
 \\
\!\!\!\! \bm \gamma^{(m)}\odot \bm d&=&\overline Q_d^TX\left( \overline \eta^{(m)}  \odot \overline \varphi \right) +\underline Q_d^TX\left( \underline \eta^{(m)} \odot F \right)\!, \quad \label{eq:reduced_mass_flow}
\end{eqnarray} 
where $\odot$ is the Hadamard product and
\begin{align}
 \!\!   F \!&=\! f\!\left(\!\overline M_d(\bm \rho^{(1)} \!+\! \bm\rho^{(2)}), \sum_m \!\sigma_m^2 ( M_s\bm s^{(m)} \!+\! M_d\bm \rho^{(m)} )\!\right)\!. \label{eq:hadamard}
\end{align}
We suppose that regulators vary slowly so the time derivative of $\overline M_d$ is insignificant, justifying its removal from Eq. \eqref{eq:reduced_partial}.  

\edit{The interested reader is encouraged to sketch a small network such as the one in Fig. \ref{bc_conf}; label the edges, nodes, compressors, and regulators; and construct the associated weighted incidence matrix, its submatrices, and their signed correspondences.  Upon substituting these matrices into Eq.  \eqref{eq:reduced_partial} and  performing the matrix multiplications, we should observe that each component of the resulting vector equation is simply Eq. \eqref{eq:pde1_integrate} corresponding to the edge component.  Likewise, the components of Eq. \eqref{eq:hadamard} reduce to Eq. \eqref{eq:pde2_integrate} after some algebraic rearrangements.}  Multiplying both sides of Eq. \eqref{eq:reduced_partial} on the left by $\overline Q_d^TX$ and using Eq. \eqref{eq:reduced_mass_flow}, we may combine Eq. \eqref{eq:reduced_partial} and Eq. \eqref{eq:reduced_mass_flow} to form the equation $\overline Q_d^TX L\overline M_d\dot {\bm \rho}^{(m)}=[Q_d^TX(\underline \eta^{(m)} \odot F)- \bm \gamma^{(m)} \odot \bm d]$, where we have used $Q_d=(\underline Q_d+\overline Q_d)$.   By writing edge concentrations in terms of nodal concentrations, and nodal concentrations in terms of concentrations of flows into the nodes, the system in Eqs. \eqref{eq:reduced_partial}-\eqref{eq:hadamard} may be written for both $m=1$ and $m=2$ as
\begin{eqnarray}
R\dot {\bm \rho}^{(m)}&=&Q_d^TX\left( \left( |\underline Q_s| \bm \alpha^{(m)}+ |\underline Q_d| \bm \eta^{(m)}\right)\odot F \right) \nonumber \\
&& \qquad - \left( I_q\bm \beta^{(m)} + I_w\bm \eta^{(m)}\right)\odot \bm d, \quad \label{eq:partial_den_sys}
\end{eqnarray} 
where $R=\overline Q_d^TX L\overline M_d$.  The system in Eq. \eqref{eq:partial_den_sys} will be called the partial density system of ODEs.  Each row $k$ of $\overline M_d$ contains exactly one nonzero component given by $\overline M_{kj}=\overline \mu_k$ for $k\in$$_{\mapsto}j$.  Using the additional fact that $X$ and $L$ are diagonal, it can be shown that the mass matrix $R$ on the left-hand-side of Eq. \eqref{eq:partial_den_sys} is diagonal with positive diagonal components given by $\bm r_{j}=\sum_{k\in _{\mapsto}j} \chi_k\ell_k\overline \mu_k$ for $j\in \mathcal V_d$.  Therefore, the matrix $R$ may readily be inverted to obtain a nonlinear control system in the usual, although complicated, ODE form.
The initial condition in Eq. \eqref{eq:pde_ic}, sampled at the refined nodes of the network, is the time-invariant solution of the system in Eq. \eqref{eq:partial_den_sys} with $\bm d=\bm d(0)$, $\bm \alpha^{(m)}=\bm \alpha^{(m)}(0)$, and $\bm \beta^{(m)}=\bm \beta^{(m)}(0)$.  We assume that this steady-state solution is the initial condition of the partial density system.

\vspace{-1ex}
\section{Equivalent Systems} \label{sec:equivalent_sys}

The system in Eq. \eqref{eq:partial_den_sys} is expressed in terms of partial densities at non-slack nodes.  Equivalent systems expressed in terms of other variables of interest may be derived from Eq. \eqref{eq:partial_den_sys} using appropriate transformations, such as that performed in the continuous case going from Eqs. \eqref{eq:pde_1a}-\eqref{eq:pde_convect} to Eqs. \eqref{eq:pde1}-\eqref{eq:pde2}.  In fact, such transformations exist even for homogeneous gas systems.  For example, the equations of natural gas flow may be expressed in terms of pressure and velocity, in terms of density and mass flux, or in terms of their dimensionless quantities.  Define vectors $\bm \rho$, $\bm p$, $\bm \nu^{(m)}$, and $\bm E$ of nodal values for density, pressure, volumetric concentration, and energy, respectively,  at non-slack nodes by  
\begin{eqnarray}
    \bm \rho&=&\bm \rho^{(1)}+\bm \rho^{(2)}, \\
    \bm p&=&\sigma_1^2\bm \rho^{(1)}+\sigma_2^2\bm \rho^{(2)}, \\
    \bm \nu^{(m)}&=&\frac{\sigma_m^2\bm \rho^{(m)}}{\sigma_1^2\bm \rho^{(1)}+\sigma_2^2\bm \rho^{(2)}}, \\
     \bm E&=&(|\overline Q_d^T| X\underline \varphi)\odot \left(r^{(1)}\bm \eta^{(1)}+r^{(2)}\bm \eta^{(2)}\right),
\end{eqnarray}
where $r^{(1)}=44.2$ (MJ/kg) and $r^{(2)}=141.8$ (MJ/kg).  Equivalent systems may be expressed in terms of any two vector variables from the set $\{\bm \rho^{(m)},\bm \eta^{(m)},\bm \nu^{(m)},\bm \rho, \bm p, \bm E \}$, excluding pairs from the subset $\{\bm \eta^{(1)},\bm \eta^{(2)},\bm \nu^{(1)},\bm \nu^{(2)}\}$ because variables in the latter subset would reduce to constant vectors in the case of homogeneous mixtures. The choice of which equivalent system to use may depend on the sought application, although some systems have better conditioning with fewer nonlinear operations than others.  Define (potentially time-varying) localized wave speed vectors $\bm a=(\sigma_1^2\bm \alpha^{(1)}+\alpha_2^2\bm \alpha^{(2)})^{1/2}$ and $\bm b=(\sigma_1^2\bm \beta^{(1)}+\sigma_2^2\bm \beta^{(2)})^{1/2}$, where the square-root is applied component-wise. The transformation from partial densities to \emph{total density and pressure} is obtained by superimposing Eq. \eqref{eq:partial_den_sys} for $m=1,2$ to obtain an equation for $\dot{\bm \rho}$ and linearly combining Eq. \eqref{eq:partial_den_sys} for $m=1,2$ with coefficients $\sigma_1^2$ and $\sigma_2^2$ to obtain an equation for $\dot{\bm p}$.  This transformation produces the system
\begin{eqnarray}
R\dot {\bm \rho}&=&Q_d^T XF - \bm d, \label{eq:pressure_density_sys1} \\
R\dot {\bm p}&=&Q_d^TX \left(\left(|\underline Q_s| \bm a^2 + |\underline Q_d| \frac{\bm p}{\bm \rho} \right) \odot F \right) \nonumber \\
&& \qquad \quad -\left(I_q \bm b^2 + I_w \frac{\bm p}{\bm \rho} \right)\odot \bm d, \quad \label{eq:pressure_density_sys2}
\end{eqnarray} 
where $F=f(\overline M_d\bm \rho, M_s \bm p_s+M_d\bm p)$.  The system in  Eqs. \eqref{eq:pressure_density_sys1}-\eqref{eq:pressure_density_sys2} will be called the total density and pressure system of ODEs.  We do not derive other equivalent systems.  Instead, we compute the solution of the partial density system of ODEs numerically, and, thereafter, obtain the other variables of interest by subsequently applying the appropriate transformations.

If $\bm \eta^{(m)}$ is a constant vector, then the system of total density and pressure decouples into two isolated subsystems that are equivalent to one another because, for constant concentration, we have $\bm p=\bm c^2 \odot \bm \rho$, where $\bm c=(\sigma_1^2\bm \eta^{(1)}+\sigma_2^2\bm \eta^{(2)})^{1/2}$ is a constant vector. In particular, the total density and pressure system in 
 Eqs. \eqref{eq:pressure_density_sys1}-\eqref{eq:pressure_density_sys2} reduces by half its dimension to the isolated system
\begin{eqnarray}
R\dot {\bm p}&=&Q_d^T X\left( \left( |\underline Q_s| \bm a^2 + |\underline Q_d| \bm c^2 \right) \right. \nonumber  \\ 
&&\qquad \quad \odot \left. f\left(\overline M_d\frac{\bm p}{\bm c^2}, M_s \bm p_s+ M_d\bm p  \right) \right) \nonumber  \\
&& -\left(I_q \bm b^2 + I_w\bm c^2 \right)\odot \bm d.  \quad \label{eq:pressure_isolated_sys}
\end{eqnarray}
The system in Eq. \eqref{eq:pressure_isolated_sys} is called the isolated total pressure system of ODEs.  Equivalent isolated subsystems expressed in terms of one vector variable from the set $\{\bm \rho^{(m)},\bm \rho, \bm E \}$ may be derived.  Each isolated subsystem is only applicable if the concentration vector $\bm \eta^{(m)}$ is constant.   Rigorous definitions and proofs of conditions on $\bm \alpha^{(m)}$, $\bm \beta^{(m)}$, $\bm q$, $\bm w$, and network topology that would result in $\bm \eta^{(m)}$ being constant are outside the scope of this study.

\section{Monotonicity} \label{sec:monotonicity}

The monotonicity properties of solutions to initial boundary value problems for flows of a homogeneous gas through an actuated transport network were examined as a means to reduce the complexity of optimization and optimal control of natural gas networks in the presence of uncertainty \cite{misra2020monotonicity}.  Here, we examine how such concepts can be extended to the transport of inhomogeneous gas mixtures, and specifically to characterize the acceptable extent and variability of hydrogen blending into a natural gas pipeline.  We first present some analytical results before proceeding with numerical simulations in the next section.

A nonlinear input-to-state initial-value system of ODEs may be generally expressed as
\begin{equation}
    \dot x = g(x,u,d), \qquad x(0)=y, \label{eq:general_non_sys}
\end{equation}
where $x(t) \in \mathcal X\subset \mathbb{R}^n$ is the state vector, $u(t)\in \mathcal U\subset \mathbb{R}^m$ is the control input vector, and $d(t)\in \mathcal D\subset \mathbb{R}^r$ is the parameter input vector defined for $t\in [0,T]$.  It is assumed that the subsets $\mathcal X$, $\mathcal U$, and $\mathcal D$ are compact and convex and that the function $g:\mathcal X \times \mathcal U \times \mathcal D \rightarrow \mathcal X$ is Lipschitz in $\mathcal X \times \mathcal U \times \mathcal D$.

{\bf Definitions:} Suppose that two independent state solutions $\{x_1(t),x_2(t)\}\subset \mathcal X$ exist (and are thus unique because $g$ is Lipschitz) with initial conditions $\{y_1,y_2\}\subset \mathcal X$, and which correspond to the piecewise-continuous control inputs $\{u_1(t),u_2(t)\} \subset \mathcal U$ and piecewise-continuous parameter inputs $\{d_1(t),d_2(t)\} \subset \mathcal D$ for $t\in [0,T]$.  For the given set of control inputs, the system in Eq. \eqref{eq:general_non_sys} is said to be {\it monotone-ordered}  with respect to $d(t)$ if $x_1(t)\le x_{2}(t)$ for $t\in [0,T]$ whenever $y_1\le y_{2}$ and $d_1(t)\le d_{2}(t)$, where inequalities for vectors are taken componentwise.  In this case, the solution states $x_1$ and $x_2$ are said to be {\it monotone-ordered}.  For simplicity, we say that a monotone-ordered system and a set of monotone-ordered solutions are {\it monotone}, {\it monotonic}, and have the property of {\it monotonicity}.  An $n\times n$ matrix $A$ is called {\it Metzler} if all of its off-diagonal elements are non-negative, i.e. $A_{ij}\ge 0$ for all $i\not = j$.  An $n\times m$ matrix is called {\it non-negative} if all of its entries are non-negative.

{\bf Theorem 1 (Monotonicity) \cite{angeli2003monotone, hirsch2006monotone}:}
The nonlinear system in Eq. \eqref{eq:general_non_sys} is monotone if and only if the Jacobian matrices $\partial g/\partial x$ and $\partial g/\partial d$ are, respectively, Metzler and non-negative almost everywhere in $\mathcal X \times \mathcal U\times \mathcal D$.

\subsection{Homogeneous Concentration} \label{sec:monotonicity_homogeneous}

The equivalent systems described in Section \ref{sec:equivalent_sys} are first reformulated in terms of the monotone system definitions above.  In steady-state \cite{misra2020monotonicity}, the pressure $\bm p$ increases componentwise with {\it decreasing} withdrawal $\bm w\ge 0$ and with {\it increasing} injection $-\bm q\le 0$.   In reference to Eq. \eqref{eq:general_non_sys}, we assume that all non-slack nodes are injection nodes  and define the input parameter by $d= \{\bm p_s,\bm d\}=\{\bm p_s,-\bm q\}$.  

{\bf Proposition 1 (Monotonicity of Total Pressure and Density):}
Assume that i) all non-slack nodes are injection nodes; ii) gas flows only in the positive direction through each edge according to its orientation in the network graph; iii) pressure is positive in each node; and iv) Eq. \eqref{eq:eps} is satisfied.  Suppose that the concentration vector $\bm \eta^{(2)}$ is constant and that there exist two state solutions $\bm p_1$, $\bm p_2$ of the system in Eq. \eqref{eq:pressure_isolated_sys} with respective initial conditions $\bm \pi_1$, $\bm \pi_2$, slack node pressures $(\bm p_s)_1$, $(\bm p_s)_2$, and non-slack injection flows $\bm q_1$, $\bm q_2$ for a given fixed set of control inputs $\{\underline \mu, \overline \mu\}$.  Here, the vector subscripts denote the first and second solutions and not the refined nodes.   If $\bm \pi_1\le \bm \pi_2$, $(\bm p_s)_1(t) \le (\bm p_s)_2(t)$,  and $\bm q_1(t)\ge \bm q_2(t)$ componentwise for all $t\in [0,T]$, then $\bm p_1(t)\le \bm p_2(t)$.  Consequently, $\bm \rho_1(t)\le \bm \rho_2(t)$, where $\bm \rho_1$ and $\bm \rho_2$ are the total densities of the two solutions.

{\bf Proof}:  Throughout this proof, the state and input subscripts correspond to the nodes of the refined graph.  Because flow is in the positive oriented direction, it follows from Eq. \eqref{eq:pde2} that $\underline \mu_k\bm p_i(t) > \overline \mu_k \bm p_j(t)$ for all $i,j\in \mathcal V$ with $k:i\mapsto j$.   Thus, the sign and absolute value operations in Eq. \eqref{eq:flux_expression} are unnecessary.  The $j$-th state dynamics in Eq. \eqref{eq:pressure_isolated_sys} for $j\in \mathcal V_d$ may be written as
\begin{align}
 \bm r_{j}\dot{\bm p}_j  &= \sum_{k:i\mapsto j} \frac{\bm \sigma^2_{i} \chi_k\Lambda_k}{\bm c_j} \left( \overline \mu_k\bm p_j \left( \underline \mu_k\bm p_i -\overline \mu_k \bm p_j \right) \right)^{1/2} \nonumber \\ 
  & - \sum_{k:j \mapsto i} \frac{\bm c^2_{j} \chi_k \Lambda_k}{\bm c_{i}} \left( \overline \mu_k\bm p_{i}\left( \underline \mu_k\bm p_j-\overline \mu_k\bm p_i \right)  \right)^{1/2} \\ & \qquad \qquad  +\bm b_j^2\bm q_j, \nonumber
\end{align}
where $\bm p_i=(\bm p_s)_i$ and $\bm \sigma_i^2=\bm a_i^2$ if $i\in \mathcal V_s$, whereas $\bm \sigma_i^2=\bm c_i^2$ if $i\in \mathcal V_d$.
It is clear from this expanded form that the function on the right-hand-side of Eq. \eqref{eq:pressure_isolated_sys} is continuously differentiable (hence Lipschitz) in the state and input variables over the domain of positive flow and pressure.  In reference to Theorem 1, we first show that the state Jacobian matrix is Metzler, i.e., $\partial \dot{\bm p_j}/\partial \bm p_i$ is non-negative for all $i,\,j\in \mathcal V_d$ with $i\not=j$. If $i$ and $j$ are non-adjacent with $i\not=j$, then clearly $\partial \dot{\bm p_j}/\partial \bm p_i=0$. Suppose that $i$ and $j$ are adjacent with $k:j\mapsto i$.  Substituting Eqs. \eqref{eq:pde_bc1}-\eqref{eq:pde_bc2} into Eq. \eqref{eq:eps} and using the relation between pressure and partial densities, we can show that $ (\overline \mu_k \bm p_i-\underline \mu_k \bm p_j)> -\overline \mu_k \bm p_i $. Thus, the Jacobian component
\begin{eqnarray}
\frac{\partial \dot{\bm p_j}}{\partial \bm p_i} &=&  \frac{\bm c_j^2\chi_k\Lambda_k \overline \mu_k (2\overline \mu_k\bm p_i-\underline \mu_k \bm p_j) }{2\bm r_j\bm c_i  (\overline \mu_k\bm p_i( \underline \mu_k\bm p_j-\overline \mu_k\bm p_i))^{1/2}}
\end{eqnarray}
is positive.  Suppose that $i$ and $j$ are adjacent with $k:i\mapsto j$.  Then
\begin{eqnarray}
\frac{\partial \dot{\bm p_j}}{\partial \bm p_i} &=&   \frac{\bm \sigma_i^2 \chi_k\Lambda_k \underline \mu_k \overline \mu_k\bm p_j  }{2\bm r_j\bm c_j  (\overline \mu_k \bm p_j( \underline \mu_k\bm p_i-\overline \mu_k\bm p_j)) ^{1/2}} > 0.
\end{eqnarray} 
Because $j\in \mathcal V_d$ is arbitrary, it follows that the state Jacobian matrix is Metzler.  We now show that the parameter Jacobian matrix is non-negative.  The above computation can be extended to show that $\partial \dot{\bm p_j}/\partial (\bm p_s)_i$ is non-negative for $i\in \mathcal V_s$.   With respect to mass inflow parameters, the Jacobian components $ \partial \dot{ \bm p_j} / \partial \bm q_i=\bm b_j^2/\bm r_j \delta_{i,j}$ are non-negative ($\delta_{i,j}$ is the Kronecker delta).   
We conclude from Theorem 1 that the system in Eq. \eqref{eq:pressure_isolated_sys} is monotone.  Because $\bm p_j=\bm c_j^2\bm \rho_j$ for $j\in \mathcal V_d$, it follows that the isolated total density system is monotone as well. $\square$

{\bf Corollary 1 (Monotonicity of Equivalent Systems):}
Assume that the conditions hold from Proposition 1.  Then $\bm \rho_1^{(m)}(t)\le \bm \rho_2^{(m)}(t)$ componentwise for all $t\in [0,T]$, where $\bm \rho_1^{(m)}$ and $\bm \rho_2^{(m)}$ are the partial densities of the two solutions.

{\bf Proof}:  The mass fraction $\bm \eta^{(m)}$ is constant, therefore it follows from Proposition 1 that  $   \bm \rho_1^{(m)}=\bm \eta^{(m)} \odot \bm \rho_1 \le \bm \eta^{(m)} \odot \bm \rho_2 = \bm \rho_2^{(m)}$. $\square$

\subsection{Heterogeneous Concentration} \label{sec:monotonicity_heterogeneous}

{\bf Proposition 2 (Non-Monotonicity of Total Pressure and Density):}
Assume that i) all non-slack nodes are injection nodes; ii) gas flows only in the positive direction through each edge according to its orientation in the network graph; and iii) pressure and density are positive in each node.   Suppose that, for a given fixed set of control inputs $\{\underline \mu, \overline \mu\}$, there exist two state solutions $(\bm \rho,\bm p)_1$, $(\bm \rho,\bm p)_2$ of the system in Eqs. \eqref{eq:pressure_density_sys1}-\eqref{eq:pressure_density_sys2} with respective initial conditions $(\bm \varrho,\bm \pi)_1$, $(\bm \varrho,\bm \pi)_2$, slack inputs $(\bm \rho_s,\bm p_s)_1$, $(\bm \rho_s,\bm p_s)_2$, and non-slack mass inflows $\bm q_1$, $\bm q_2$ that satisfy $(\bm \varrho,\bm \pi)_1\le (\bm \varrho,\bm \pi)_2$, $(\bm \rho_s(t),\bm p_s(t))_1\le (\bm \rho_s(t),\bm p_s(t))_2$,  and $\bm q_1(t)\ge \bm q_2(t)$ componentwise for all $t\in [0,T]$.  If $\bm \eta^{(m)}(t)$ is time-varying, then, in general, $(\bm \rho(t),\bm p(t))_1 \not \le(\bm \rho(t),\bm p(t))_2$ component-wise for all $t\in [0,T]$.

 {\bf Proof}:  
Throughout this proof, the state and input subscripts correspond to the nodes of the refined graph.  From Theorem 1, it suffices to show that one component of the state Jacobian matrix  is negative.  The $j$-th nodal pressure dynamics in Eq. \eqref{eq:pressure_density_sys2} may be written as
\begin{align}
 \bm r_{j}\dot{\bm p}_j  &= \sum_{k:i\mapsto j} \bm \sigma^2_i\chi_k\Lambda_k \left( \overline \mu_k \bm \rho_j \left( \underline \mu_k \bm p_i-\overline \mu_k \bm p_j  \right)\right)^{1/2} \nonumber \\
  & - \sum_{k:j \mapsto i}\frac{\bm p_j}{\bm \rho_j} \chi_k\Lambda_k \left( \overline \mu_k\bm \rho_i \left( \underline \mu_k \bm p_j -\overline \mu_k\bm p_i\right)  \right)^{1/2} \\ & \qquad \qquad +\bm  b_j^2\bm q_j, \nonumber
\end{align}
where $\bm p_i=(\bm p_s)_i$, $\bm \rho_i=(\bm \rho_s)_i$, and $\bm \sigma_i^2=\bm a_i^2$ if $i\in \mathcal V_s$, and $\bm \sigma_i^2=\bm p_i/\bm \rho_i$ if $i\in \mathcal V_d$. By adding a refined edge to the graph if necessary, we assume that there is an edge $k':i'\mapsto j$ with $i'\in \mathcal V_d$.  The Jacobian component corresponding to $\bm \rho_{i'}$ is given by
\begin{eqnarray*}
\frac{\partial \dot{\bm p_j}}{\partial \bm \rho_{i'}} &=&  -\frac{\chi_{k'} \Lambda_{k'}}{\bm r_j}\frac{\bm p_{i'}}{\bm \rho_{i'}^2} \left( \overline \mu_{k'} \bm \rho_j \left( \underline \mu_{k'} \bm p_{i'}-\overline \mu_{k'}\bm p_j  \right)\right)^{1/2},
\end{eqnarray*}
which is negative.  It follows from Theorem 1 that the system in Eqs. \eqref{eq:pressure_density_sys1}-\eqref{eq:pressure_density_sys2} is not monotone, regardlesss of Eq. \eqref{eq:eps}. $\square$

\section{Network Case Study} \label{sec:network}

\begin{figure}
\centering
\includegraphics[width=\linewidth]{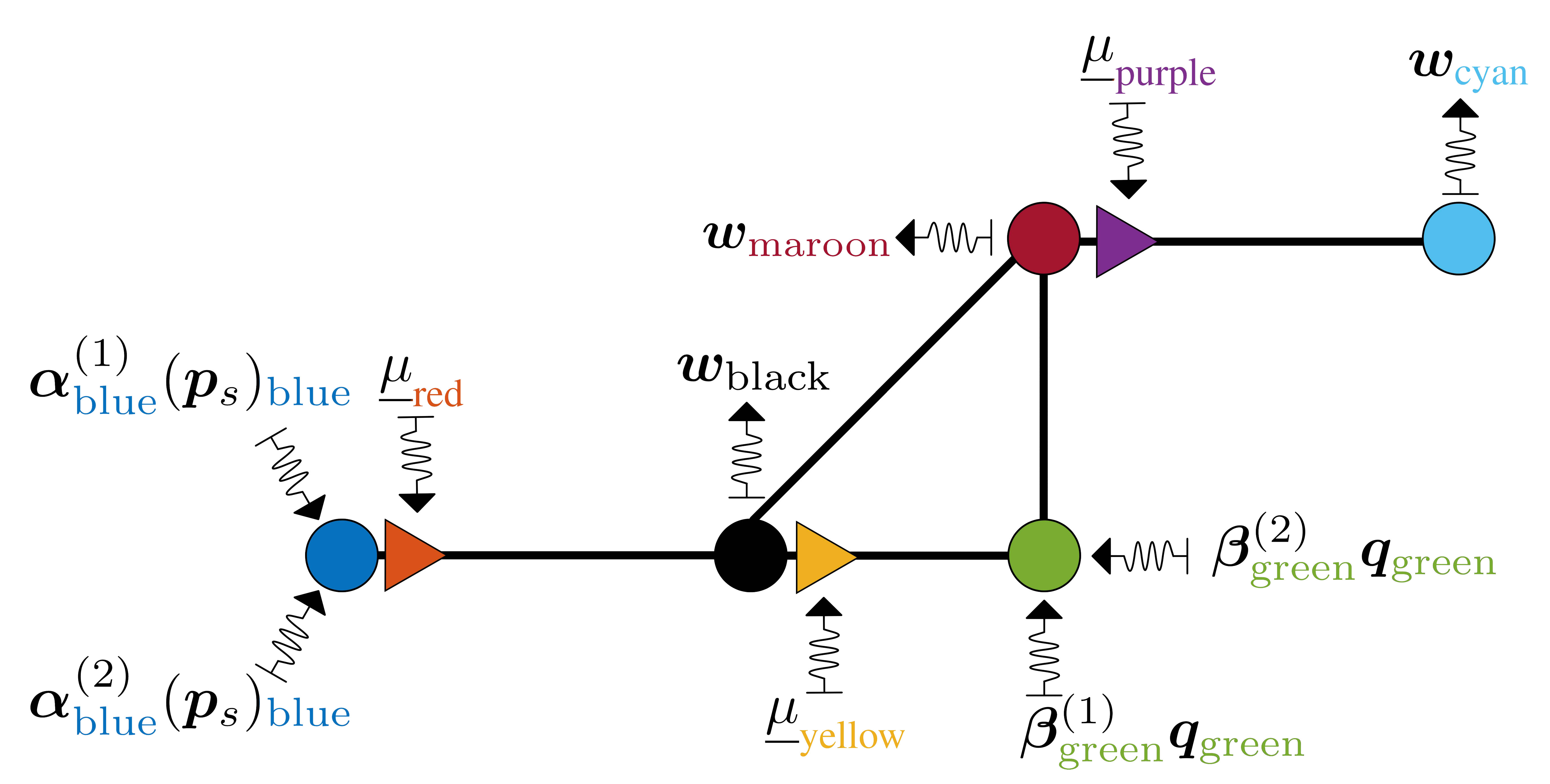}
\caption{Network configuration (not to scale).  The triangles represent compressor stations.
    Pipeline dimensions between nodes: blue to black (20 km), black to green (70 km), green to maroon (10 km), black to maroon (60 km), maroon to cyan (80 km).  The pipelines have uniform diameter (0.9144 m) and friction factor (0.01), except for the black to maroon pipeline that has diameter (0.635 m) and friction factor (0.015).}
\label{net_conf}
\end{figure}

\begin{figure}
\centering
\includegraphics[width=\linewidth]{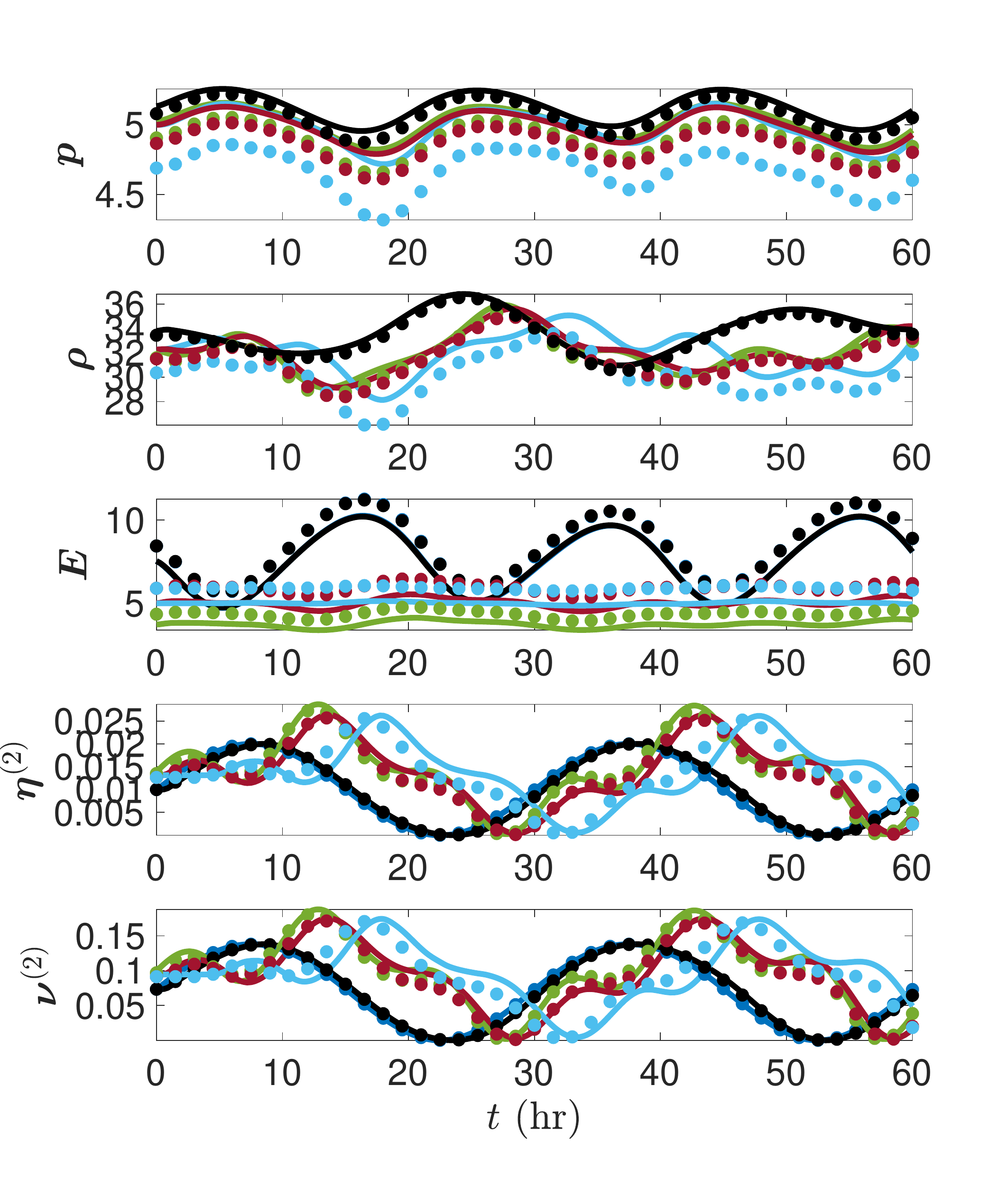}
\caption{Two solutions (solid lines vs. dots) at the color-coordinated network nodes in Fig. \ref{net_conf}.  From top to bottom, the depicted nodal solutions are total pressure (MPa), total density (kg/m$^3$), total energy (GJ/s), concentration by mass, and concentration by volume.  The boundary conditions for both solutions are $(\bm p_s)_{\text{blue}}=5$ MPa, $\bm \alpha_{\text{blue}}^{(2)}(t)=0.01(1+\sin(4\pi t/T))$, $\bm \beta_{\text{green}}^{(2)}(t)=0.125(1+\sin(12\pi t/T))$, $\bm q_{\text{green}}(t)=3$ (kg/s), $\bm w_{\text{black}}(t)=60(1-\sin(6\pi t/T))$ (kg/s),  $\mu_{\text{red}}=1.0678$, $\underline \mu_{\text{yellow}}=1.0140$, and $\underline \mu_{\text{purple}}=1.0734$, where $T=60$ (hr).  The boundary condition that differs between the two solutions is $\bm w_{\text{cyan}}(t)=110$ (kg/s) (solid lines) and $\bm w_{\text{cyan}}(t)=130$ (kg/s) (dots). } 
\label{net_ex1}
\end{figure}

\begin{figure}
\centering
\includegraphics[width=\linewidth]{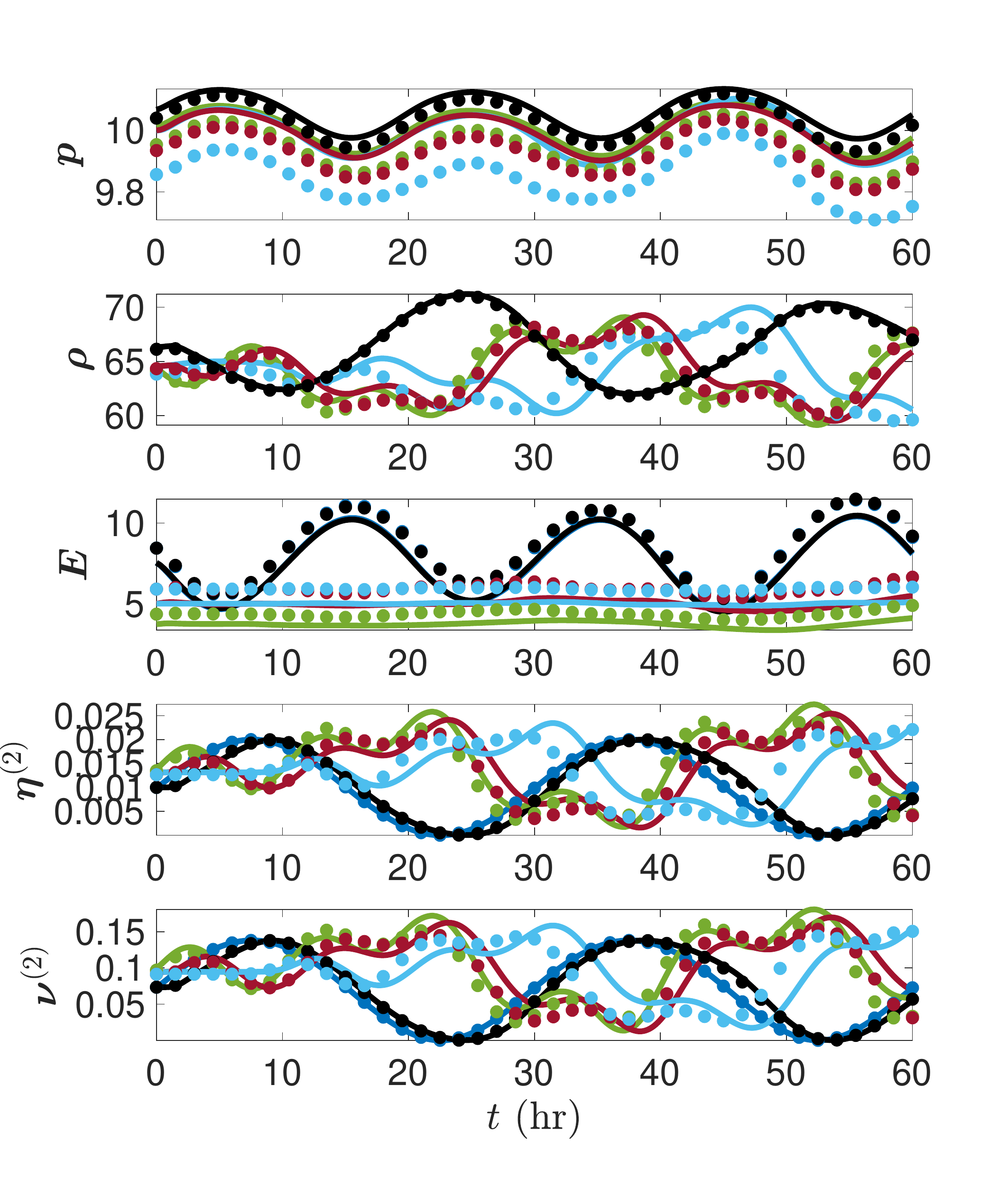}
\caption{Same boundary conditions as in Fig. \ref{net_ex1} except for $(\bm p_s)_{\text{blue}}=10$ (MPa).}
\label{net_ex2}
\end{figure}

\begin{figure}
\centering
\includegraphics[width=\linewidth]{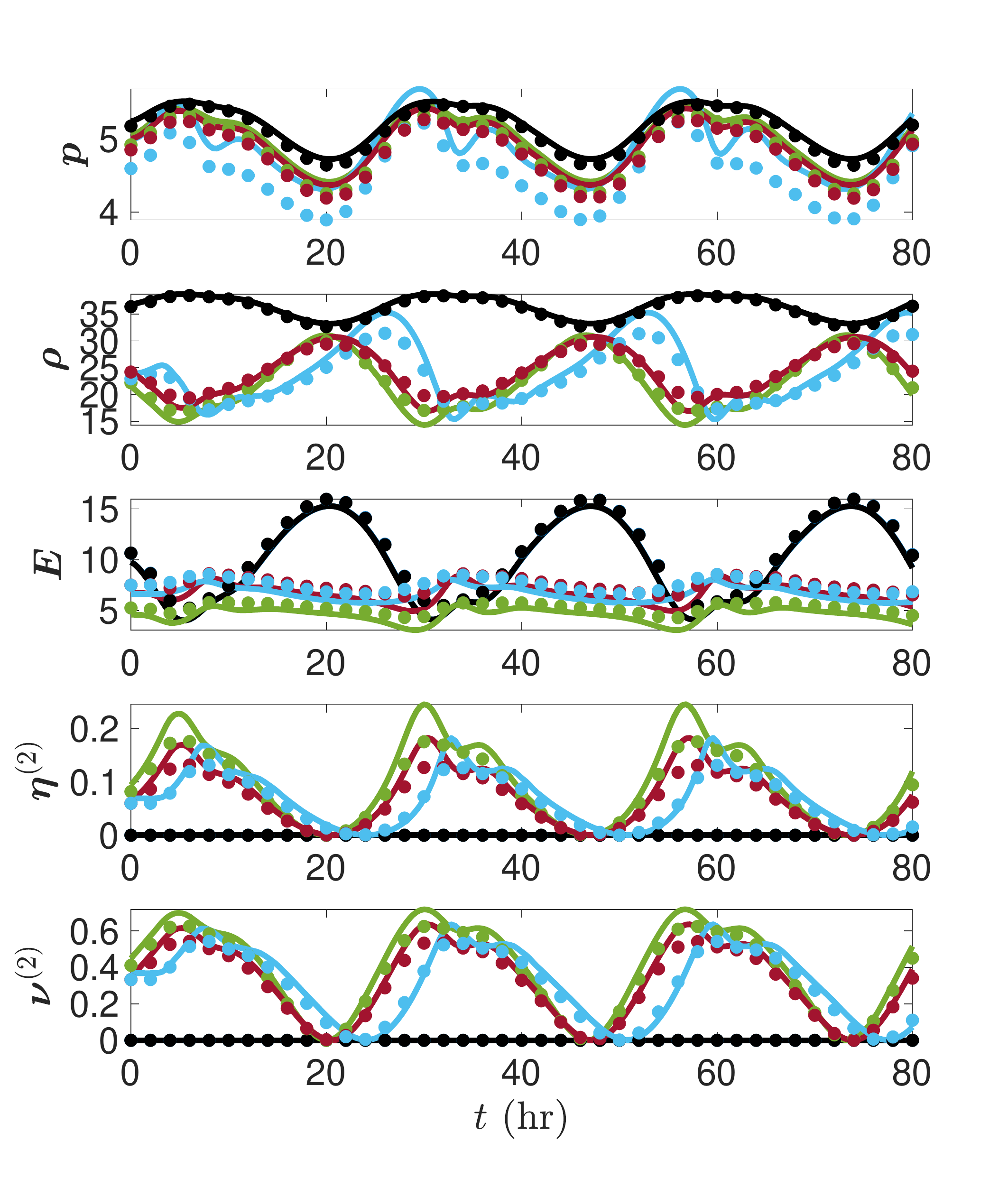}
\vspace{-4ex}
\caption{The boundary conditions for the two solutions are $(\bm p_s)_{\text{blue}}=5$ (MPa), $\bm \alpha_{\text{blue}}^{(2)}(t)=0$, $\bm \beta_{\text{green}}^{(2)}(t)=1$, $\bm q_{\text{green}}(t)=9(1+\sin(6 \pi t/T))$ (kg/s), $\bm w_{\text{black}}(t)=100(1-\sin(6\pi t/T))$ (kg/s), $\mu_{\text{red}}=1.1096$, $\underline \mu_{\text{yellow}}=1.0057$, and $\underline \mu_{\text{purple}}=1.1301$, where $T=80$ (hr).  The other boundary condition is $\bm w_{\text{cyan}}(t)=130$ (kg/s) (solid lines) and $\bm w_{\text{cyan}}(t)=150$ (kg/s) (dots).}
\label{net_ex3}
\end{figure}

\begin{figure}
\centering
\includegraphics[width=\linewidth]{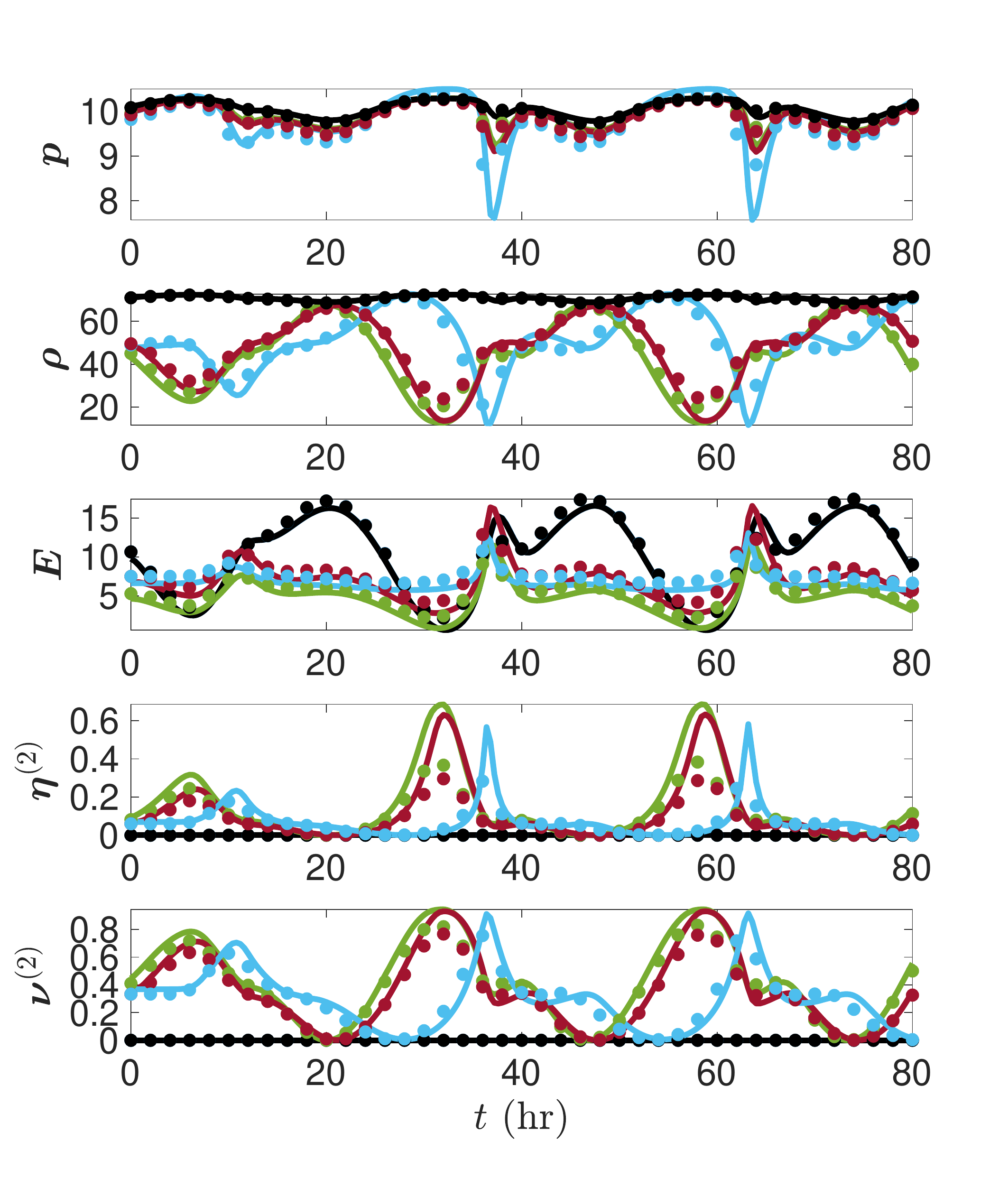}
\caption{Same boundary conditions as in Fig. \ref{net_ex3} except for $(\bm p_s)_{\text{blue}}=10$ (MPa).} 
\label{net_ex4}
\end{figure}

\begin{figure}
\centering
\includegraphics[width=\linewidth]{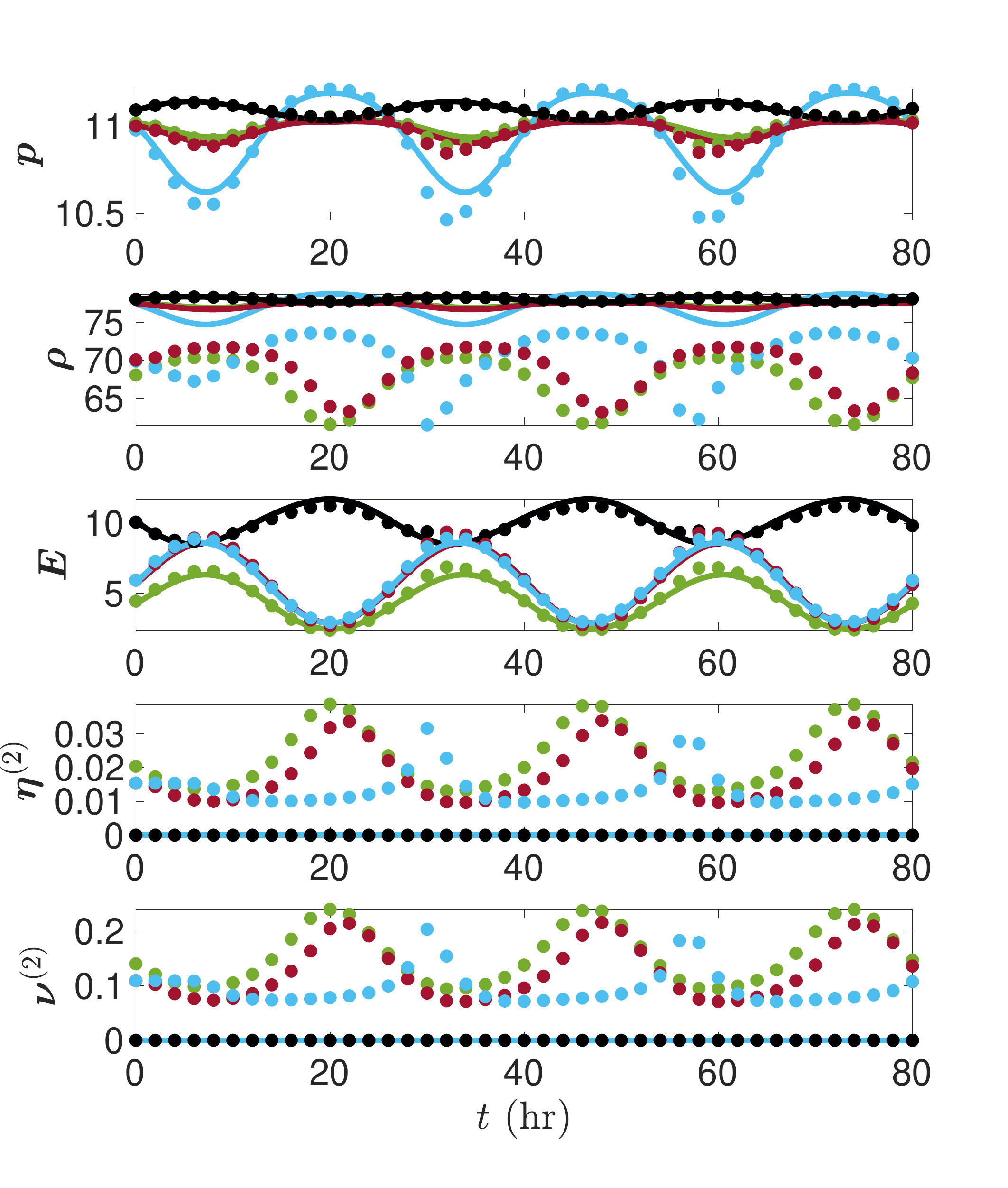}
\caption{The boundary conditions for the two solutions are $(\bm p_s)_{\text{blue}}=11$ (MPa), $\bm \alpha_{\text{blue}}^{(2)}(t)=0$, $\bm \beta_{\text{green}}^{(2)}(t)=1$, $\bm w_{\text{black}}(t)=100(1-\sin(6\pi t/T))$ (kg/s), $\bm w_{\text{cyan}}(t)=130(1+0.5\sin(6 \pi t/T))$ (kg/s), $\mu_{\text{red}}=1.0240$, $\underline \mu_{\text{yellow}}=1.0029$, and $\underline \mu_{\text{purple}}=1.0199$, where $T=80$ (hr).  The other boundary condition is $\bm q_{\text{green}}(t)=0$ (kg/s) (solid lines) and $\bm q_{\text{green}}(t)=2$ (kg/s) (dots).}
\label{net_ex5}
\end{figure}

We use numerical simulations to examine the effects of time-varying heterogeneity of a transported mixture on the flow dynamics of equivalent system variables throughout a network.  The simulations are performed for a test network that was used in a previous study \cite{GYRYA201934}, in which the authors presented a staggered grid discretization method for the numerical solution of homogeneous natural gas pipeline flow.  We refer the reader to Appendix \ref{sec:comparison} in which we show the results of our implementation of the IBVP posed in the former study in order to verify that we obtain the same solution when no hydrogen is present.  The configuration and dimensions of the network are shown in Fig. \ref{net_conf}.  The dark blue node is a slack node at which pressure and concentration are specified, the black, maroon, and cyan nodes are non-slack withdrawal nodes, and the green node is a non-slack injection node.  The wave speeds are chosen to be $\sigma_1=377$ (m/s) and $\sigma_2=2.8 \sigma_1$.  We simulate several examples  to illustrate that some physical quantities may exhibit fewer crossings than others in certain response regimes, given ordered boundary parameters.  These examples provide insight into which equivalent system may be more useful for simulation and optimization of gas mixture dynamics \cite{baker2023optimal} for a particular response regime.  Figs. \ref{net_ex1}-\ref{net_ex5} show the solutions of five different examples.  Each example computes two solutions that correspond to slightly different but monotone-ordered boundary conditions.  For each example, the two solutions are depicted in the figures distinctly with solid lines and dots and with colors coordinated to match the colors of the nodes of the network diagram in Fig. \ref{net_conf}. We now describe the simulation results for each example.

The first example in Fig. \ref{net_ex1} considers time-varying sinusoidal forcing in the concentration of the gas mixtures being injected into the slack node (blue node) and the non-slack injection node (green node).  In this figure, the total pressure, density, and energy solutions at the non-slack nodes do not cross, but the mass and volumetric concentrations do show crossings.  The solutions in Fig. \ref{net_ex2} have the same boundary conditions as those in Fig. \ref{net_ex1} except for the slack node pressure.  If the slack node pressure is doubled, then the total density trajectories will exhibit crossings at each non-slack node but the pressure and energy trajectories still do not cross.  In Figs. \ref{net_ex3} and \ref{net_ex4}, the blue node injects pure natural gas and the green node injects pure hydrogen with a varying mass inflow profile.  As seen in Fig. \ref{net_ex3}, the pressure and energy solutions at each node do not show crossings.  However, a close examination shows that the density solutions exhibit crossovers at every node upstream from the point of hydrogen injection.  Moreover, the concentration solutions show crossings at only the cyan node.  The solutions in Fig. \ref{net_ex4} have the same boundary conditions as those in Fig. \ref{net_ex3} except for the slack node pressure.   If the slack node pressure is doubled, then the resulting pressure, density, and energy trajectories will cross at all of the non-slack nodes.  Moreover, the concentration trajectories in this example cross at every node upstream from the node of hydrogen injection.  At nodes downstream the injection of hydrogen, the concentration of hydrogen is zero, as it ought to be.  We note that the solutions in Figs. \ref{net_ex3} and \ref{net_ex4} may not be realistic in the current operation of natural gas pipelines because the concentration of hydrogen reaches very high levels.  However, these figures indicate that flow dynamics may exhibit similarly rapid transients in pipelines that are upgraded to deliver significant amounts of hydrogen.  All of the solution variables show large gradient surges in small time intervals.

From the simulations presented thus far, it may appear that the concentration variables are the most likely of the equivalent system variables to violate monotonicity.  The simulation shown in Fig. \ref{net_ex5} demonstrates a counterexample to this conjecture in which pressure, density, and energy trajectories all exhibit crossings even though the concentration solutions do not.  However, strictly speaking, density trajectories cross only upstream from the point of hydrogen injection.  The difference between the solid line and dotted solutions in Fig. \ref{net_ex5} is that the solid line represents the solution for homogeneous natural gas and the dotted solution represents a slightly perturbed solution that results from a small injection of hydrogen made at the green non-slack node. The concentration variables associated with the two solutions may be identical at times in certain network nodes, but there cannot be strict crossings in the concentration trajectories for this example.  The reason is because the homogeneous natural gas solution corresponds to zero hydrogen concentration and this is the lower bound that the concentration variables can achieve.

In each of the examples above, the five edges of the network are discretized into 240 refined edges with $\ell_k=1$ (km) for all $k\in \hat{\mathcal E}$ and the simulations are performed using the partial density system of ODEs in Eq. \eqref{eq:partial_den_sys}.  The equivalent system variables are computed using the transformations presented in Section \ref{sec:equivalent_sys}.  Although one kilometer discretization size is sufficient to demonstrate non-monotonicity for slowly-varying concentrations, a much smaller discretization size is required to accurately simulate rapidly-varying concentrations.  We note that even the slowly-varying solutions in Figs. \ref{net_ex1}-\ref{net_ex5} show noticeable convergence as the discretization size is decreased from 1 (km) to 100 (m).  For small discretization lengths ($\ell_k\le 100$ (m)), the crossings of the solutions in Figs. \ref{net_ex1}-\ref{net_ex5} may be more pronounced.  We will return to this point later on in Sections \ref{sec:PI} and \ref{sec:chaotic} when we simulate pipeline flows of highly heterogeneous gas mixtures.
 


\section{Monotonic Interface} \label{sec:MI}

Proposition 2 shows that the total pressure and density system of ODEs is not monotone-ordered over the entire boundary condition parameter region $\mathcal D=\{\bm p_s,\bm d,\bm \alpha^{(m)},\bm \beta^{(m)}\}$.  However, by Proposition 1 and the continuity of solutions with respect to initial conditions and inputs \cite{khalil2002nonlinear}, the non-isolated total pressure and density system of ODEs is  expected to be monotone-ordered over a certain subregion of $\mathcal D$ that consists of concentration vectors that are uniformly close to a constant concentration vector.  Again by continuity, monotonicity is also expected to hold for slow variations in concentration with large amplitudes.  This suggests that there may be a nontrivial monotonic interface (MI) that partitions $\mathcal D$ into monotonic and non-monotonic phase regions for each equivalent system variable.   Moreover, the simulation results from Section \ref{sec:network} suggest that the MIs for each equivalent system variable may be significantly different from one another.  We focus on partitioning the subregion of $\mathcal D$ that consists of only concentration boundary parameters, because the concentration variable is the only factor that leads to conditions in which monotonicity does not hold.

We analyze the MI numerically for a single pipeline with concentration and pressure specified at the inlet of the pipeline (node 1) and with mass outflow specified at the outlet (node 2).  The parameters and boundary conditions that do not change are pipeline length $\ell=50$ (km), diameter $D=0.5$ (m), friction factor $\lambda=0.11$, and constant slack node pressure $\bm p_s=7$ (MPa).  We denote the concentration of hydrogen at the inlet slack node by $\bm \alpha_1(t)=\bm \alpha_1^{(2)}(t)$ and specify it to be
\begin{equation}
    \bm \alpha_1(t)= \alpha_1 \left( 1+\kappa\sin(2\pi \omega_* t) \right), \label{eq:sine_concentration}
\end{equation}
 where $\kappa$ is the amplitude factor of the sinusoid, $\omega_*$ is the forcing frequency in cycles per hour, and $ \alpha_1$ is the mean concentration profile around which the sinusoid oscillates.  Here, the subscript is with respect to the node number.  The sub-region of $\mathcal D$ that we consider consists of all pairs $(\omega_*,\kappa)$ with $0\le \omega_* \le 2$ and $0\le \kappa \le 1$. \edit{We note that the pipeline system in Eq. \eqref{eq:partial_den_sys} with the boundary conditions specified as above may be written as an autonomous system by extending the state space by two dimensions and writing $\bm \alpha_1$ as a state solution of the harmonic oscillator.  We note this extended formulation so that our subsequent results can be interpreted in the context of autonomous dynamical systems theory. However, the extension is not necessary for the analysis, so we omit the details.}



\begin{figure}
\centering
\includegraphics[width=\linewidth]{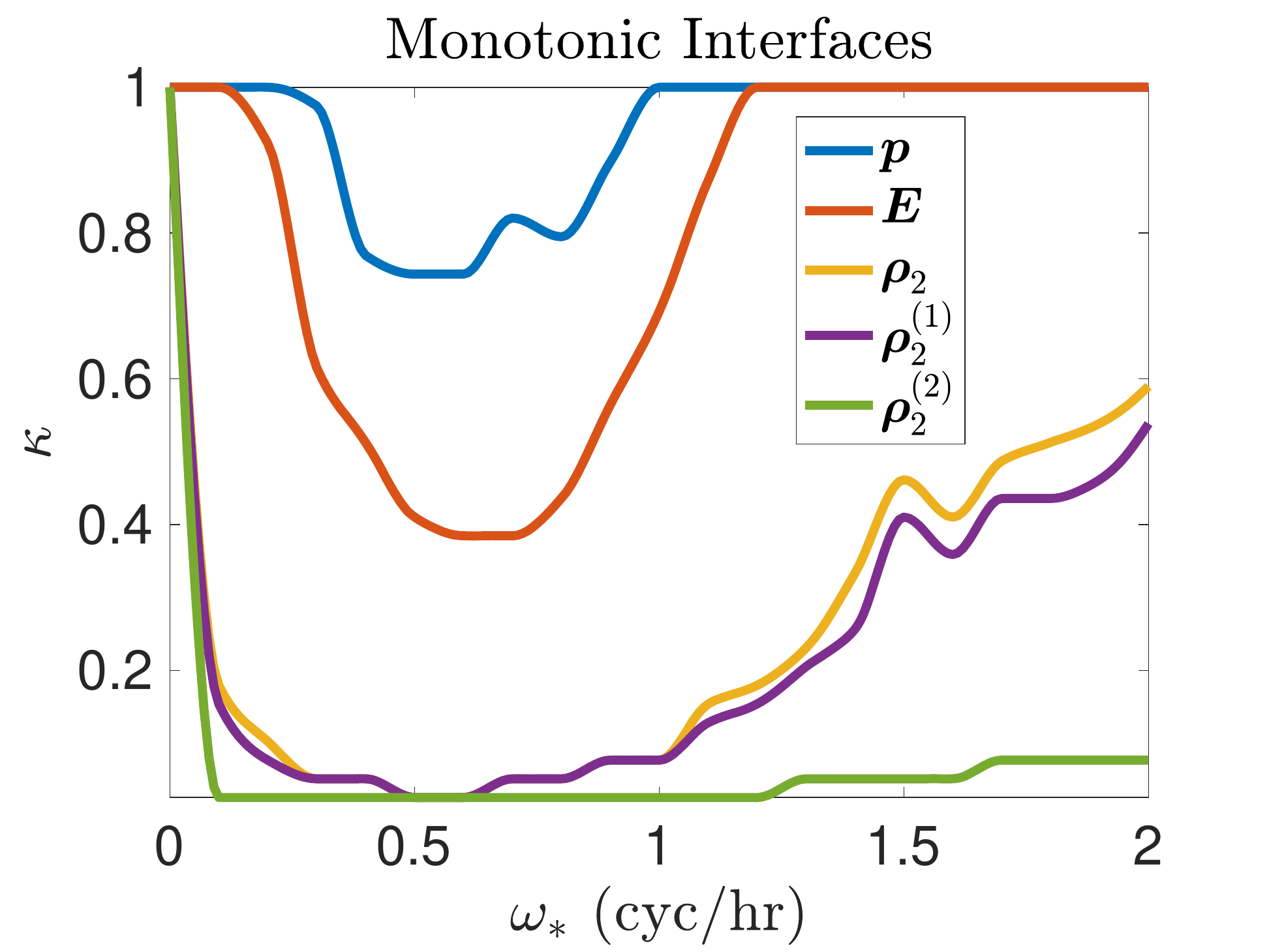}
\caption{\edit{Monotonic interfaces $(\omega_*,\kappa_*(\omega_*))$ in the $(\omega_*,\kappa)$ plane for each of the equivalent system variables listed in the legend. Two solutions corresponding to two ordered mass outflows and the concentration pair $(\omega_*,\kappa)$ will not cross one another if $\kappa<\kappa_*(\omega_*)$.}}
\label{fig:monotonic_interface}
\end{figure}

In addition to the boundary conditions and parameters specified above, we use the parameter values $\sigma_1=377$ (m/s), $\sigma_2=2.8 \sigma_1$, and $\alpha_1=0.02$ in the analyses in this section.   We now describe our process of computing the MI.  For each $(\omega_*,\kappa)$ in Eq. \eqref{eq:sine_concentration}, we compute three solutions corresponding to three monotone-ordered constant mass outflows $\bm w_2=\overline \varphi \pi(D/2)^2$ (kg/s), where $\overline \varphi=120$, 140, and 160 (kg/m$^2$s).  The region in the $(\omega_*,\kappa)$ plane defined by $0\le \omega_*\le 2$ and $0\le \kappa \le 1$ is discretized into a $21\times 41$ grid of discrete pairs.  For each discrete $\omega_*$, we compute the three solutions for each discrete $\kappa$ starting from $\kappa=0$ and increasing $\kappa$ until we achieve the lower bound $\kappa=\kappa_*(\omega_*)$ with which at least two of the three solutions exhibit crossings at some point in time.  Thus, the three solutions corresponding to the pair $(\omega_*,\kappa)$ will not cross if $\kappa<\kappa_*(\omega_*)$ and at least two of the three will cross if $\kappa\ge \kappa_*(\omega_*)$.  The MI is defined by the set of pairs $(\omega_*,\kappa_*(\omega_*))$ in the plane, and each of the equivalent system variables has its own MI associated with it.   The cubic spline interpolated MI curves for several equivalent system variables are depicted in Fig. \ref{fig:monotonic_interface}.  The region below the MI curve is called the monotone response region (MRR).  Fig. \ref{fig:monotonic_interface} shows that the MRRs for hydrogen density, natural gas density, total density, energy, and pressure are nested increasing sets where the MRRs for hydrogen density and pressure are the smallest and largest sets, respectively.  

For time-varying concentration profiles, Fig. \ref{fig:monotonic_interface} suggests that the pressure and energy equivalent system should be used if monotonicity properties are important to the formulation.  This is the conclusion that we arrive at in Section \ref{sec:network}.  Of the five examples from that section, the only examples that consider similar sinusoidal forcing in concentration are those that correspond to Figs. \ref{net_ex1} and \ref{net_ex2}.  Both of these examples introduce a combination of two sinusoidal forcing frequencies to the network, $0.1$ (cyc/hr) and 0.033 (cyc/hr).  Recall from Figs. \ref{net_ex1} and \ref{net_ex2} that only the pressure and energy solutions did not exhibit crossings for these two examples.  This observation agrees with the MIs in Fig. \ref{fig:monotonic_interface}, where the operating points $(\omega_*,\kappa)=(0.1,1)$ and $(\omega_*,\kappa)=(0.033,1)$ are above all of the MIs except for the pressure and energy MIs.  It is important to mention that the MIs computed here for a single pipeline do not necessarily match the MIs that would result for the network topology from Section \ref{sec:network}, even if the slack node pressure, total outflow, and other boundary condition parameters of the network were equal to those used for the single pipeline.

\edit{Observe that the MRRs for the pressure and energy variables are significantly larger in area than the MRRs for density, as seen in Fig. \ref{fig:monotonic_interface}.  This is inherent because of the manner in which the boundary conditions have been specified.  Recall that the boundary condition at the slack node of the pipeline presented here has been specified to maintain constant pressure.  Although the model has been derived more generally, constant slack node pressure is a common specification for pipeline simulation, so we have followed this convention.  However, if the pressure is constant and the concentration varies at the slack node, then the density will be forced to vary at the slack node, because of the equation of state of the mixture.  It follows that the density will oscillate if the concentration oscillates, and there will be some propagation of density oscillations at speeds related to the mass outflows. Because we specify different outflows (that are monotone-ordered) to compute the MIs, the density waves corresponding to the different outflows will typically be out of phase and the solutions will exhibit crossings. Therefore, we expect the density variables to be the most sensitive to ordering properties under such boundary conditions, and this is apparent in Fig. \ref{fig:monotonic_interface}.  
}

As $\omega_*$ increases from $\omega_*=0$ to $\omega_*=2$ (cyc/hr), the MI curves qualitatively decrease from unity to a lower bound, flatten out, and then increase.  The fact that the amplitude factor generally increases along the MI as $\omega_*$ increases beyond $\omega_*=0.75$ is a robustness feature of monotonicity to high frequency uncertainty.  This property appears to be a consequence of wave attenuation in strongly dissipative gas pipeline flow \cite{baker2021analysis}.  In particular, the gas pipeline demonstrates low-pass filtering characteristics with which the amplitudes of high frequency travelling waves tend to be significantly attenuated over short distances, and, therefore, the likelihood of conditions in which monotonicity does not hold decreases as the frequency of the high frequency oscillation increases.   If the concentration of hydrogen injected into the network contains a small variation of high frequency uncertainty, then the MIs demonstrate that this uncertainty typically will not cause an otherwise theoretically monotonic operation to become non-monotonic.

\section{Periodic Interface}  \label{sec:PI}

We demonstrate that non-periodic solutions can arise from sinusoidal forcing in concentration.  To numerically study periodicity and the breakdown thereof, we must simulate the solutions over long time intervals that span hundreds of hours.  In addition, we will consider large and fast variations in concentration.  As we have mentioned at the end of Section \ref{sec:network}, fast variations require an extremely fine spatial discretization size for the finite volume discretization method presented in Section \ref{sec:discretize_flow}.  The small spatial discretization size creates a large ODE system, which is difficult to implement over a long time interval on a digital computer.  Therefore in our study of periodicity, instead of using the finite volume method, we discretize the pipeline at the (translated) nodes of Chebyshev polynomials for which exponential convergence properties are obtained (e.g., see \cite{ascher2011first}).  We briefly outline the method in Appendix \ref{sec:chebyshev}.  The analysis presented in this section is performed in the single 50 (km) pipeline that was used in the previous section to study the MI, with $D=0.5$ (m), $\lambda=0.11$, and constant slack node pressure $\bm p_s=7$ (MPa).  However, in addition to the other parameters specified in that section, we now use $\sigma_1=338.38$ (m/s), $\sigma_2=4\sigma_1$, and $\alpha_1=0.2$, but all of the other parameters remain the same.

\begin{figure}
\centering
\includegraphics[width=\linewidth]{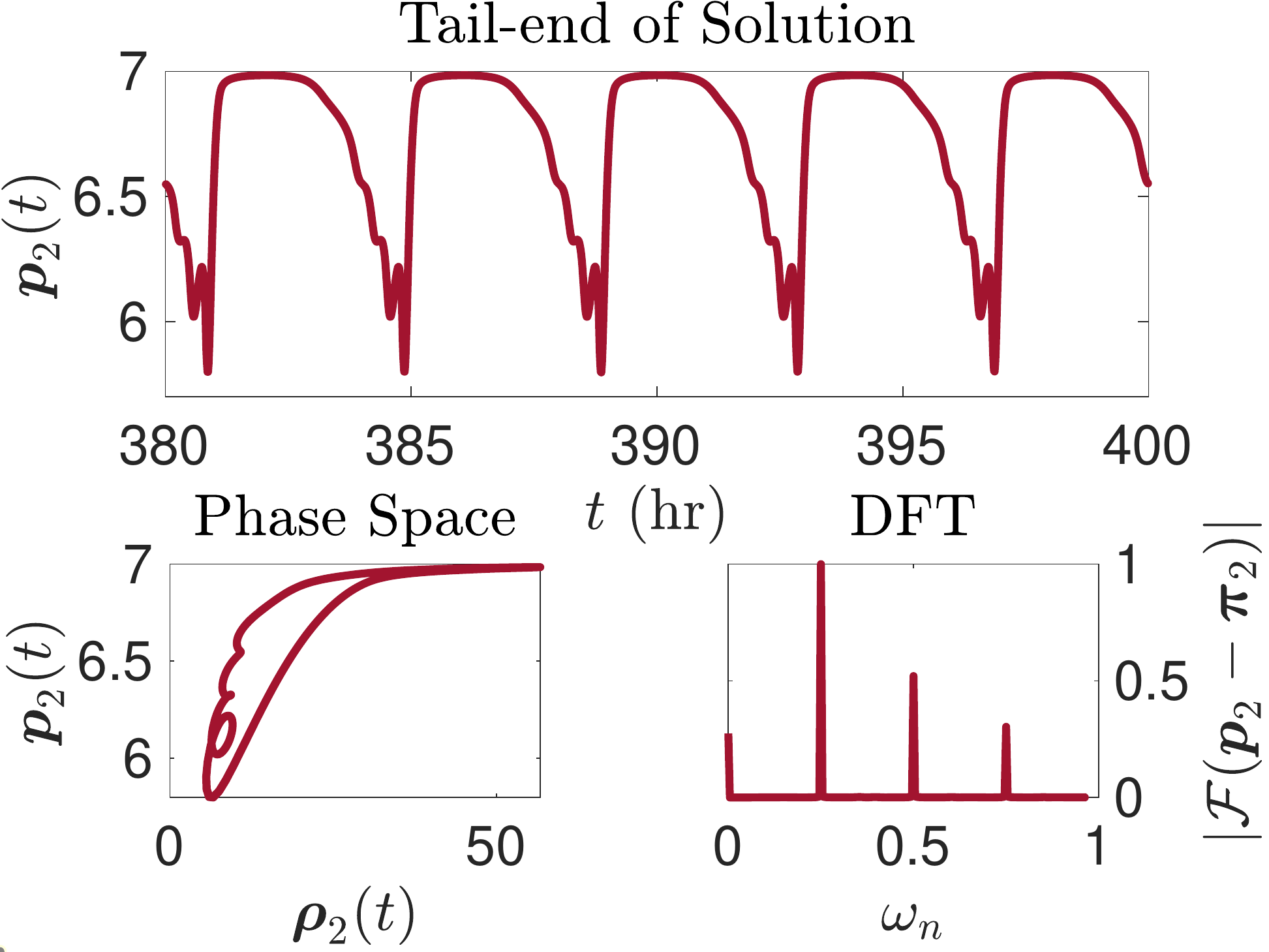}
\caption{\edit{Pipeline solution with boundary conditions $\bm w_2(t)=75 \pi(D/2)^2$ (kg/s),  $\omega_*=0.25$, and $\kappa=1.0$. The periodicity measure in Eq. \eqref{eq:power_spectrum} is $\mathcal P =0.30$.}}
\label{fig:phase2}
\end{figure}

\begin{figure}
\centering
\includegraphics[width=\linewidth]{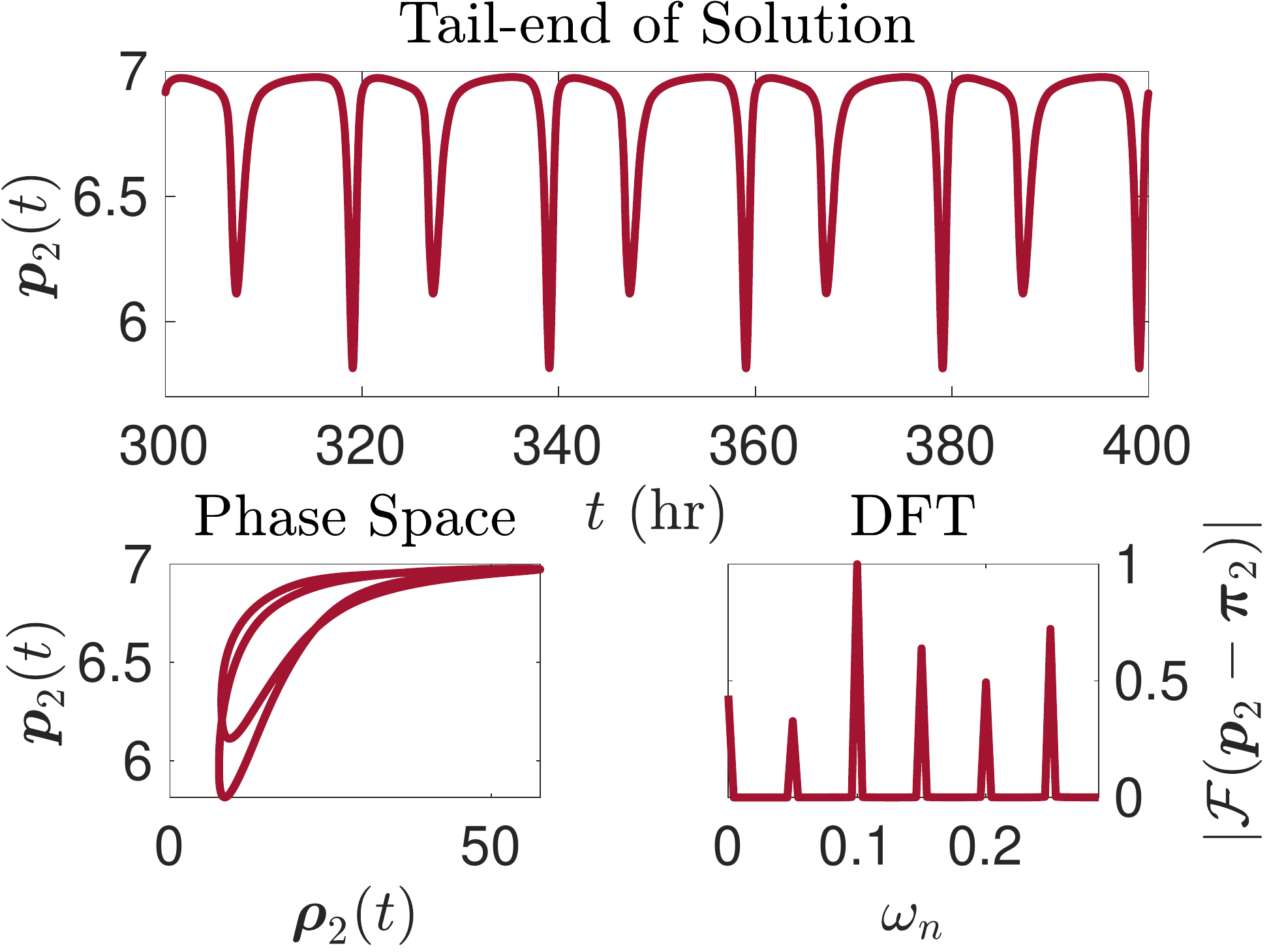}
\caption{\edit{Pipeline solution with boundary conditions $\bm w_2(t)=75 \pi(D/2)^2$ (kg/s),   $\omega_*=0.1$, and $\kappa=0.98$.  The periodicity measure in Eq. \eqref{eq:power_spectrum} is $\mathcal P=0.62$.}}
\label{fig:phase3}
\end{figure}

\begin{figure}
\centering
\includegraphics[width=\linewidth]{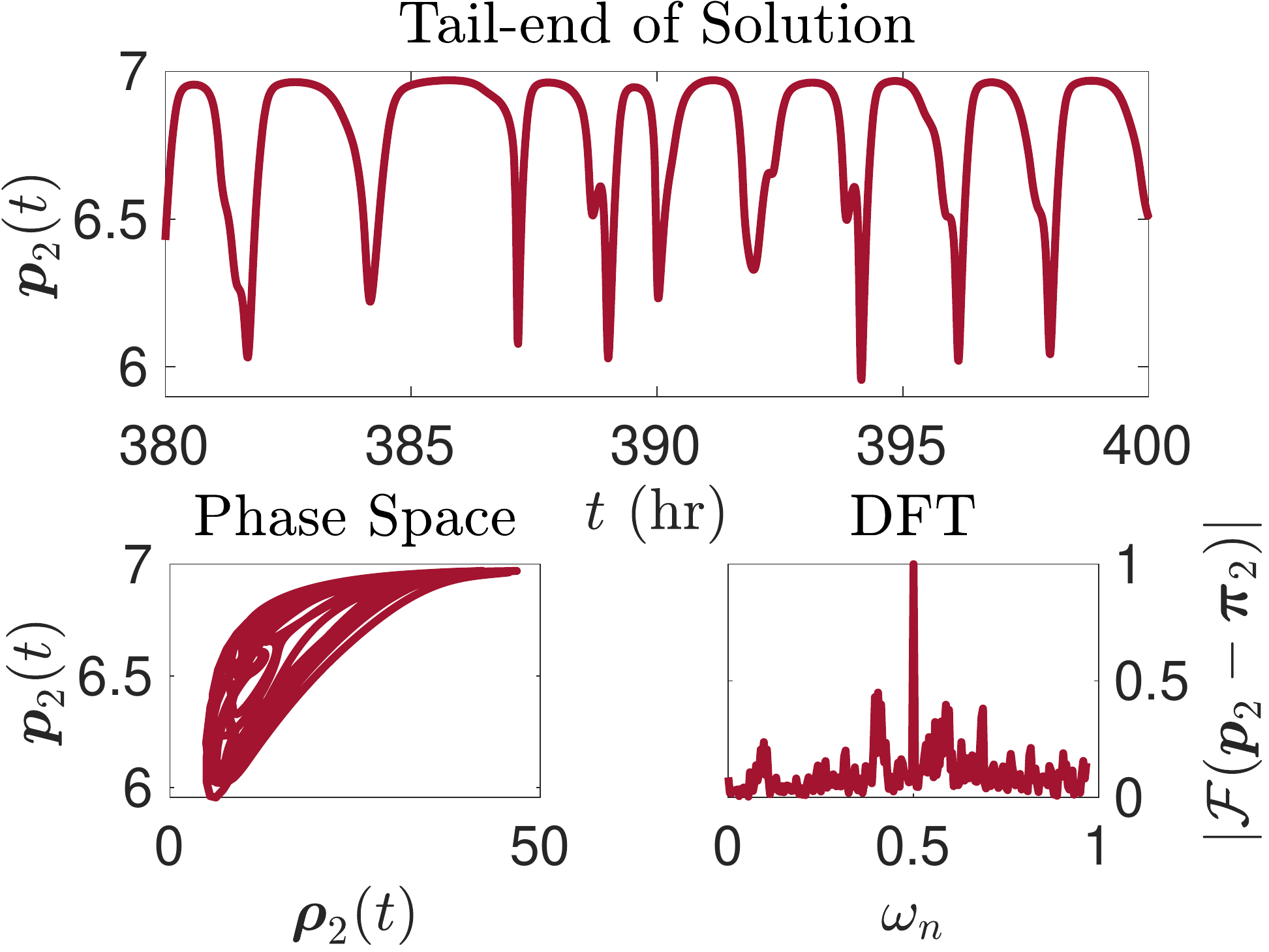}
\caption{\edit{Pipeline solution with boundary conditions $\bm w_2(t)=75 \pi(D/2)^2$ (kg/s),   $\omega_*=0.5$, and $\kappa=0.9$.  The periodicity measure in Eq. \eqref{eq:power_spectrum} is $\mathcal P=1.19$.}}
\label{fig:phase4}
\end{figure}

To introduce our analysis on periodicity, we show three examples in Figs. \ref{fig:phase2} to \ref{fig:phase4} that share the same boundary conditions with one another except for the different frequencies $\omega_*$ and amplitude factors $\kappa$ of the sinusoidal concentration profile in Eq. \eqref{eq:sine_concentration}.  The top of the three figures depict the pressure solutions at the outlet of the pipeline for $t \in [rT,T]$ with $0.7\le r\le 0.95$, where $T=400$ (hr).  The tail-ends of the solutions are used so that initial transient responses do not affect the analysis of periodic orbits.  The bottom left-hand-sides of Figs. \ref{fig:phase2} to \ref{fig:phase4} show the phase space diagrams of outlet density and outlet pressure during the later stages of the simulations.  We see that the solutions in Figs. \ref{fig:phase2} and \ref{fig:phase3} approach periodic orbits and that the solution in Fig. \ref{fig:phase4} does not appear to do so.  However, even the two periodic responses in Figs. \ref{fig:phase2} to \ref{fig:phase3} have certain properties that are not observed in homogeneous natural gas simulations \cite{baker2021analysis}.  Particularly, in Fig. \ref{fig:phase2}, multiple local minima in the pressure appear for every local minimum of the sinusoidal forcing over the time interval $[0.95T,T]$.  The additional local minima correspond to the inner loop of the periodic orbit.  The pressure in Fig. \ref{fig:phase3} has the same number of local minima as the forcing sinusoid over the interval $[0.75T,T]$, but has twice the period.   These examples demonstrate that periodic solutions arising from heterogeneous mixtures of gases may be irregular in the following sense.  From the laws of fluid dynamics, gas pressure should decrease with decreasing density under constant temperature and volume.  However, because of the oscillating gas composition, the phase space diagram in Fig. \ref{fig:phase2} contains sequences of four small time intervals during which density decreases while pressure increases, and the phase space diagram in Fig. \ref{fig:phase3} contains two such time intervals.  \edit{Note that the phase space diagram of outlet pressure and density for flow of a single ideal gas in the transient regime is simply a line with positive slope, which the solution traverses.}

The frequency responses of the outlet pressures are depicted on the bottom right-hand-sides of Figs. \ref{fig:phase2} to \ref{fig:phase4} using the discrete Fourier transform \cite{oppenheim2001discrete} defined below in Eq. \eqref{dft}.  The dominant frequency mode in the solution appears at the forcing frequency $\omega_n=\omega_*$ in Figs. \ref{fig:phase2} to \ref{fig:phase4}.  The generated harmonic modes in Fig. \ref{fig:phase2} appear at integer multiples of $\omega_*$.  This behavior is typical for homogeneous natural gas pipeline flow \cite{baker2021analysis}. The generated harmonic modes in Fig. \ref{fig:phase3} appear at half the values of the integer multiples of $\omega_*$.  This behavior indicates period-doubling bifurcations \cite{zhang2009route} at the forcing frequency $\omega_*=0.1$ as the amplitude factor $\kappa$ increases.   The pressure dynamics in Fig. \ref{fig:phase4} appear to be comprised by a continuous distribution of generated harmonic modes.  These observations inspire a quantitative measure of periodicity in terms of the frequency response of the solution.  This is the approach taken in \cite{tziperman1994nino} for the transition to chaotic responses in oceanic wind bursts.  We define a sequence of evenly-spaced samples of the tail-end of the outlet pressure by $\bm p_2[k]=\bm p_2((0.6+k/N)T)$ for $k= 0,\dots,0.4N$, where $N$ is equal to the number of time samples of the numerical solution over the interval $[0,T]$.  For such a sampled sequence $\bm \psi[k]$, the normalized \edit{discrete Fourier transform (DFT) is defined as}
\begin{equation}
    \{\mathcal F \bm \psi\}[\omega_n] = \frac{\sum_{k=0}^{0.4N}  \bm \psi[k]e^{-\bm j2\pi \omega_n k  } }{\max_{\omega_n} \left|\sum_{k=0}^{0.4N}\bm  \psi[k]e^{-\bm j2\pi \omega_n k}\right|}, \label{dft}
\end{equation}
where $\bm j$ is the imaginary unit and $\omega_n=n/(0.4T)$ (cyc/hr) are the sampling frequencies for $n=0,\dots,0.4N$.  Periodicity is measured with the average power spectrum defined by
\begin{eqnarray}
    \mathcal P = \frac{1}{0.4N+1}\sum_{n=0}^{0.4N} \left| \{\mathcal F (\bm p_2-\bm \pi_2)\}[\omega_n] \right|^2\times 100, \label{eq:power_spectrum}
\end{eqnarray}
where $\bm \pi_2=\bm p_2(0)$ is the initial steady-state value of pressure at the outlet of the pipeline.   The shifted pressure in the power spectrum is used to suppress the zero frequency component of the initial state.

\begin{figure}
\centering
\includegraphics[width=\linewidth]{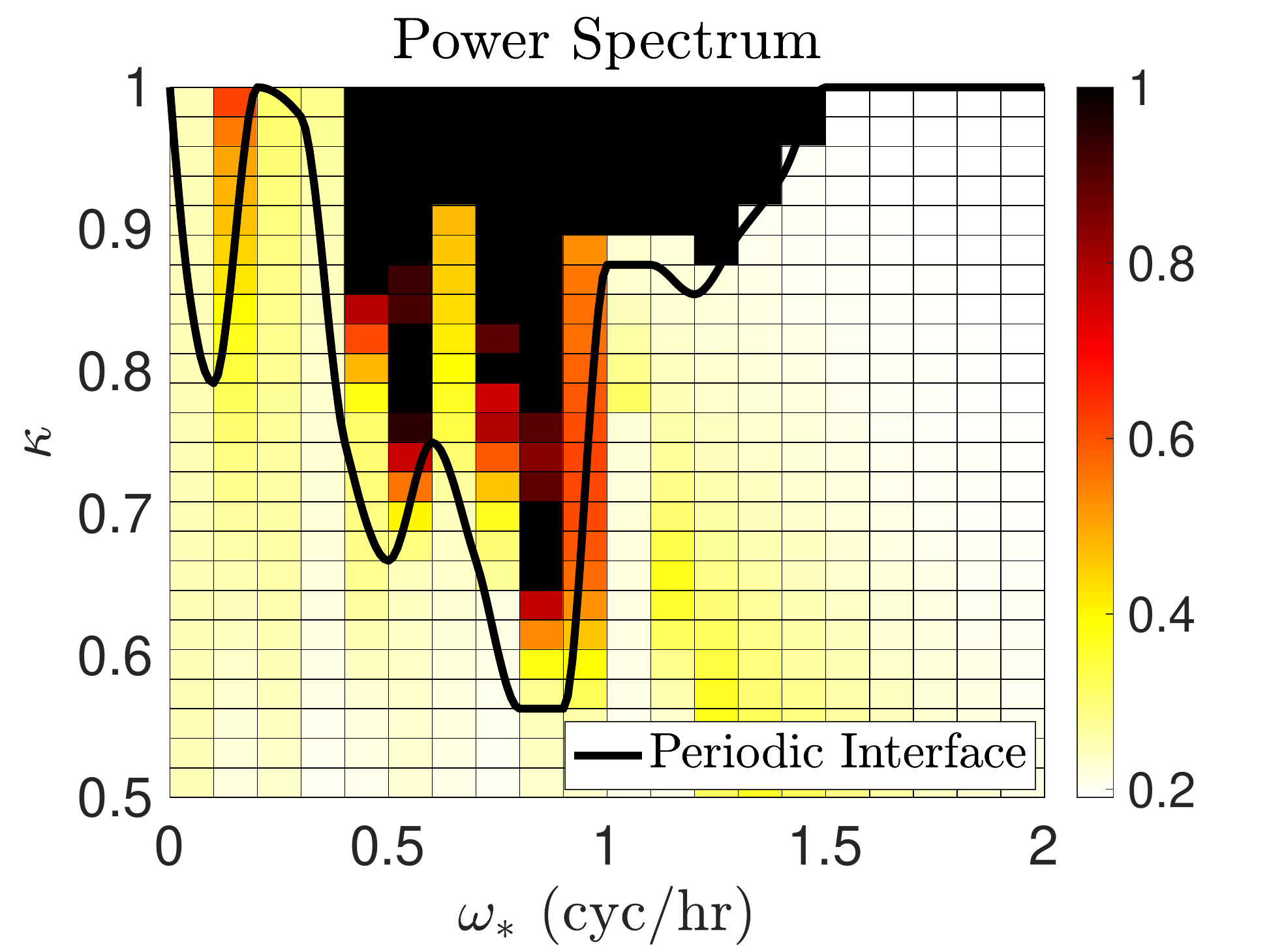}
\caption{Color map of the power spectrum $\mathcal P$ in Eq. \eqref{eq:power_spectrum} as a function of $(\omega_*,\kappa)$ in Eq. \eqref{eq:sine_concentration}. The boundary conditions are $ \alpha_1=0.2$ and $\bm w_2(t)=75$.   In this figure, we plot the minimum between 1 and $\mathcal P$ in Eq. \eqref{eq:power_spectrum}.}
\label{fig:powerspectrum}
\end{figure}

The power spectrum $\mathcal P$ is depicted in the form of a color map as a function of $(\omega_*,\kappa)$ in Fig. \ref{fig:powerspectrum}, where $\omega_*$ is the forcing frequency and $\kappa$ is its amplitude factor given in Eq. \eqref{eq:sine_concentration}.   This figure has been obtained numerically as follows.  Similarly to the way that we have computed the MIs, the region in the $(\omega_*,\kappa)$ plane defined by $0\le \omega_*\le 2$ and $0.5\le \kappa \le 1$ is discretized into a $21\times 25$ grid of discrete pairs.  For each frequency and amplitude factor of the forcing concentration on this grid, we numerically simulate the solution in the pipeline for 400 hours.  We then compute the normalized DFT and the average power spectrum of the tail-end of the sampled solution as defined above.  These computations provide the discrete set of quantified values depicted in Fig. \ref{fig:powerspectrum}.   \edit{The periodic interface (PI) in Fig. \ref{fig:powerspectrum} partitions the $(\omega_*,\kappa)$ plane into periodic and non-periodic response subregions.  The computation of the PI is performed as follows.  For each $\omega_*$, the parameter $\kappa$ is increased from $\kappa=0$ to $\kappa=\kappa^*(\omega_*)$, where  $\kappa^*(\omega_*)$ is the upper bound on $\kappa$ below which the DFT of the outlet pressure consists of countably many pulses.   From numerical simulations, the average power spectrum corresponding to the forcing pair $(\omega_*,\kappa^*(\omega_*))$ is typically around $\mathcal P=0.3$ (which depends on parameters and the number of time samples of the numerical solution). Therefore, we quantitatively define $\kappa^*(\omega_*)$ to be the upper bound on $\kappa$ below which $\mathcal P<0.3$. The periodic response region (PRR) is defined to be the set of boundary condition parameter pairs $(\omega_*,\kappa)$ with $\kappa<\kappa^*(\omega_*)$.}

\vspace{-1ex}
\section{Chaotic Interface} \label{sec:chaotic}

\begin{figure}
\centering
\includegraphics[width=\linewidth]{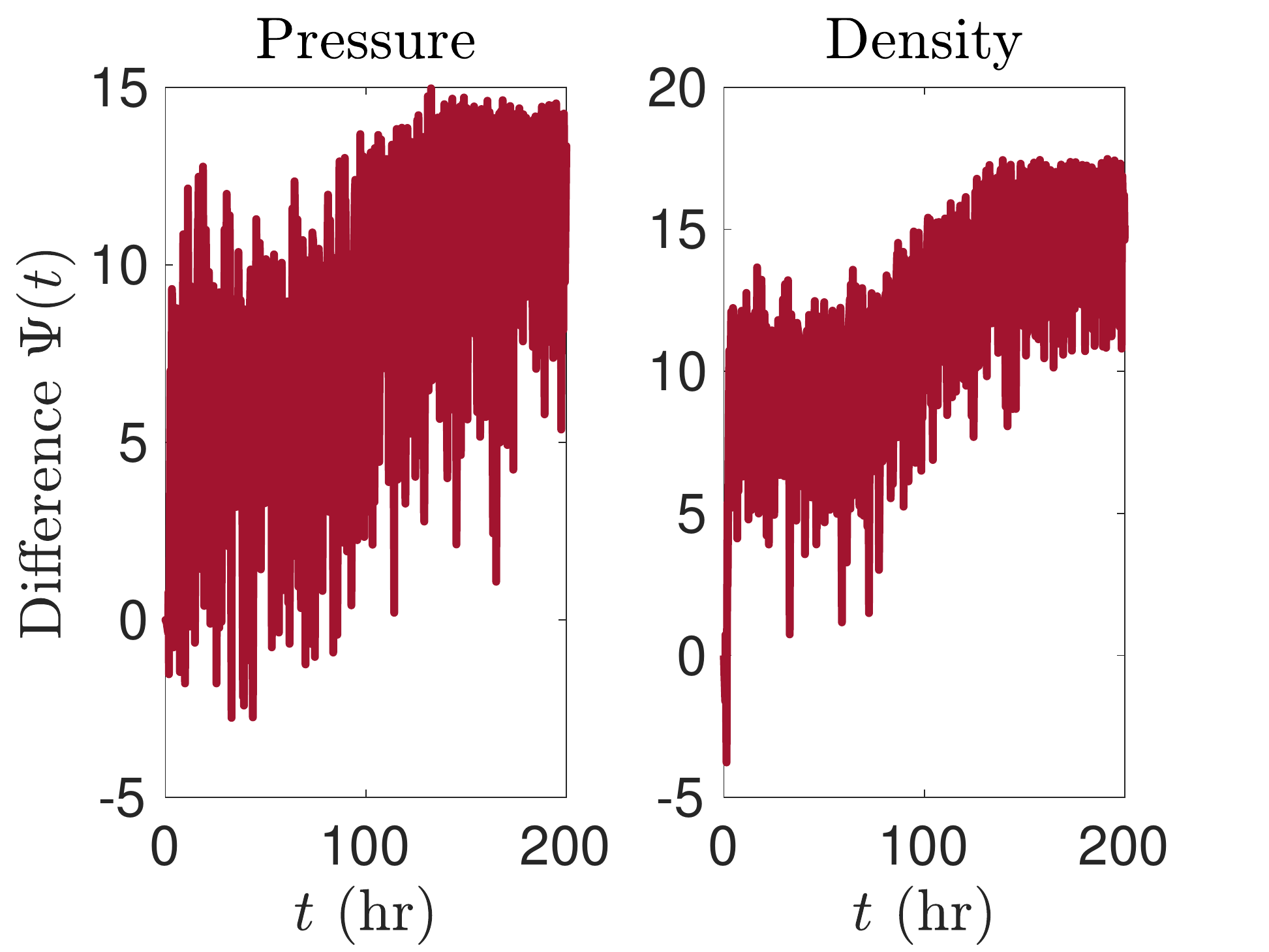}
\caption{\edit{Example of the difference $\Psi[n]=\Psi(t_n)$ between two outlet pressure solutions and corresponding outlet density solutions in the pipeline from Section \ref{sec:PI}. The boundary condition parameters are $(\omega_*,\kappa)=(0.5,0.95)$, $\alpha_1=0.2$, and the initial conditions of the two solutions correspond to steady mass flows $\bm w_2(0)=75(\pi D^2/4)$ and $\bm w_2(0)=75.0001(\pi D^2/4)$.}}
\label{fig:difference}
\end{figure}

\edit{The analysis presented in the previous section demonstrates that heterogeneous mixtures of gases may behave irregularly for periodic boundary conditions with frequency and amplitude $(\omega_*,\kappa)$ outside of the periodic response region.  In this section, we demonstrate that the solutions that exhibit such behavior are also chaotic in the sense of being highly sensitive to initial conditions \cite{lorenz1995essence}.  The extent of chaos in a finite-dimensional system can be quantified by the largest Lyapunov exponent of the system \cite{benettin1980lyapunov}.  This measure provides an estimate on the exponential rate of divergence between two solutions that begin their trajectories close to one another.  Several approaches have been developed to approximate the Lyapunov exponents of an unknown finite-dimensional system from numerical data by using the time series associated with a numerical solution of the system \cite{wolf1985determining, rosenstein1993practical, brown1991computing}.  These methods require an appropriate embedding dimension on which the calculation of the Lyapunov exponent depends.  Although there are methods that approximate the optimal embedding dimension of time-series observations of a finite-dimensional system  \cite{kennel1992determining}, there may not always be a consistent approximation for an observable generated by an underlying infinite-dimensional PDE system.  We measure the extent of chaos in our system by using a rate of divergence between two specific solutions.  It is important to note that the measure we present does not necessarily provide an estimate of the largest Lyapunov exponent of our system.

Our analysis will be performed in the single pipeline with parameters and boundary conditions defined in the previous section.  Consider two scalar-valued time series solutions $\psi_1[n]$ and $\psi_2[n]$ with $|\psi_2[0]-\psi_1[0]|<\delta$, where $\delta$ is small relative to $\psi_1[0]$.   Here, $\psi[n]=\psi(t_n)$ represents any one of the dynamic variables evaluated at the outlet of the pipe at discrete times $t_n=(n/N)T$ for $n=0,\dots,N$, where $(N+1)$ is the number of time samples of the numerical solution over the time interval $[0,T]$.  The rate of divergence between the two solutions is defined to be the slope of a linear estimation of the difference function
\begin{equation*}
    \Psi[n]= \log\left| \frac{\psi_2[n]-\psi_1[n]}{\psi_2[0]-\psi_1[0]} \right|
\end{equation*}
over an interval $I_0=[n_0,n_1]$ during which the two solutions diverge exponentially (assuming that they diverge).  Fig. \ref{fig:difference} shows an example of how the difference function may increase with time for two divergent solutions that began their trajectories close to one another.  Of course, boundedness of the physical system limits the divergence to within a bounded region in the state space.  For the example in Fig. \ref{fig:difference}, the rate of divergence may be computed as the slope of a linear estimation of $\Psi[n]$ over the interval $[0.5N,0.55N]$ corresponding to $t_n\in [100,110]$.  One such linear estimation is the least squares regression line.

The interval $I_0$ typically varies with respect to parameters and could be difficult to construct for two general solutions that may or may not diverge.  For two solutions that do not diverge, the slope of the linear estimation of $\Psi[n]$ on the interval $[0,N]$ would be at most approximately zero. Moreover, there are usually several subintervals $I_0\subset [0,N]$ over which the slope of the linear estimation of $\Psi$ changes sign, particularly for two solutions that do not diverge. Consequently, the linear estimation of $\Psi[n]$ over an interval $I_0$ may not always reflect the actual mean growth of $|\psi_2[n]-\psi_1[n]|$ over the entire interval $[0,N]$.  In an attempt to moderate these difficulties, the amount of divergence between two solutions $\psi_1$ and $\psi_2$ is defined by the measure
\begin{eqnarray} \label{eq:chaos}
    \mathcal C :&=& \frac{1}{n_2-n_1}\left(\frac{1}{|I_T|}\sum_{n\in I_T}\Psi[n] -  \frac{1}{|I_0|}\sum_{n\in I_0}\Psi[n] \right) \\
                &=&\frac{1}{(n_2-n_1)|I_T|}\sum_{n\in I_T}\log|\psi_2[n]-\psi_1[n]| \nonumber \\
                && \qquad   -  \frac{1}{(n_2-n_1)|I_0|}\sum_{n\in I_0}\log|\psi_2[n]-\psi_1[n]|,  \qquad \nonumber
\end{eqnarray}
with $I_T=[n_2,n_3]$ where $n_1<n_2$.  The function $\mathcal C$ in Eq. \eqref{eq:chaos} will be referred to as the chaos measure.  This measure is interpreted as the slope of the local mean of $\Psi$ or as the slope of the local mean logarithmic difference of the two solutions computed over an initial interval $I_0$ and a final interval $I_T$.  Large and positive $\mathcal C$ indicates exponential divergence between the two time series (over some interval).  The means over the initial and final time intervals are used to suppress fluctuations and estimate the expected values of the associated quantities over the intervals.  The interval $I_0$ must span a range that precedes any exponential divergence and should also not include initial transients.  Moreover, the interval $I_T$ must follow any exponential divergence.

\begin{figure}
\centering
\includegraphics[width=\linewidth]{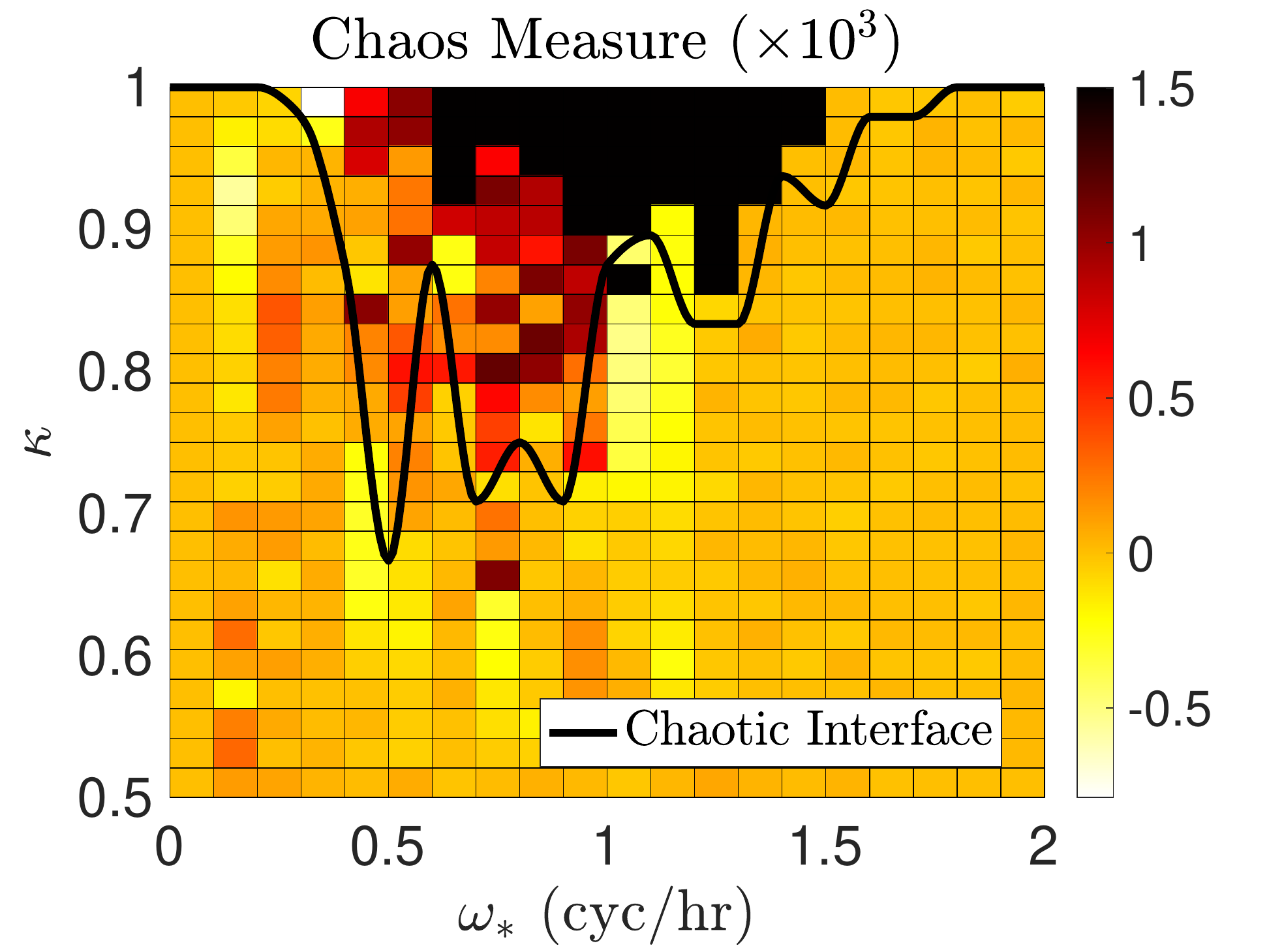}
\caption{\edit{Color map of the measure $\mathcal C$ in Eq. \eqref{eq:chaos} as a function of $(\omega_*,\kappa)$ in Eq. \eqref{eq:sine_concentration}. The boundary conditions are the same as those in Fig. \ref{fig:powerspectrum} and the initial conditions of the two solutions correspond to $\bm w_2(0)=75(\pi D^2/4)$ and $\bm w_2(0)=75.1(\pi D^2/4)$.}}
\label{fig:chaos_measure}
\end{figure}

We analyze the measure $\mathcal C$ for the outlet pressure $\psi[n]= p(t_n,\ell)$ of the single pipeline used previously in Section \ref{sec:PI} with $\alpha_1=0.2$, but we note that equivalent system variables share similar divergence properties.    As for the computations of the power spectrum shown in Fig. \ref{fig:powerspectrum}, we discretize the region in the $(\omega_*,\kappa)$ plane defined by $0\le \omega_*\le 2$ and $0.5\le \kappa \le 1$ into a $21\times 25$ grid of discrete pairs.  For each frequency and amplitude factor of the forcing concentration on this grid, we numerically simulate two solutions in the pipeline for 100 hours whose initial conditions correspond to withdrawal rates $\bm w_2(0)=75(\pi D^2/4)$ and $\bm w_2(0)=75.1(\pi D^2/4)$, where the subscript refers to the outlet node.   We then compute the value $\mathcal C$ for each discrete pair on the grid using the intervals $I_0=[0.08N,0.15N]$ and $I_T=[0.5N,0.8N]$ with $N=10,000$ and then depict the collected values as a color map in Fig. \ref{fig:chaos_measure}.   While we suppose that the intervals $I_0$ and $I_T$ are sufficient for the specified initial conditions, these intervals may not be adequate for the example in Fig. \ref{fig:difference} because the solutions in that example are initially closer to one another and require a larger time to begin diverging. For each frequency, the values of $\mathcal C$ in Fig. \ref{fig:chaos_measure} are more scattered than the associated values of $\mathcal P$ in Fig. \ref{fig:powerspectrum}.  In the latter case, the values of the power spectrum either increase or decrease almost monotonically for each frequency as the values of the amplitude factor increase.  However, for a fixed frequency, the values of the chaos measure may oscillate sporadically around $\mathcal C=0$ as $\kappa$ increases from $\kappa=0$ to $\kappa=1$. 

\begin{figure}
\centering
\includegraphics[width=\linewidth]{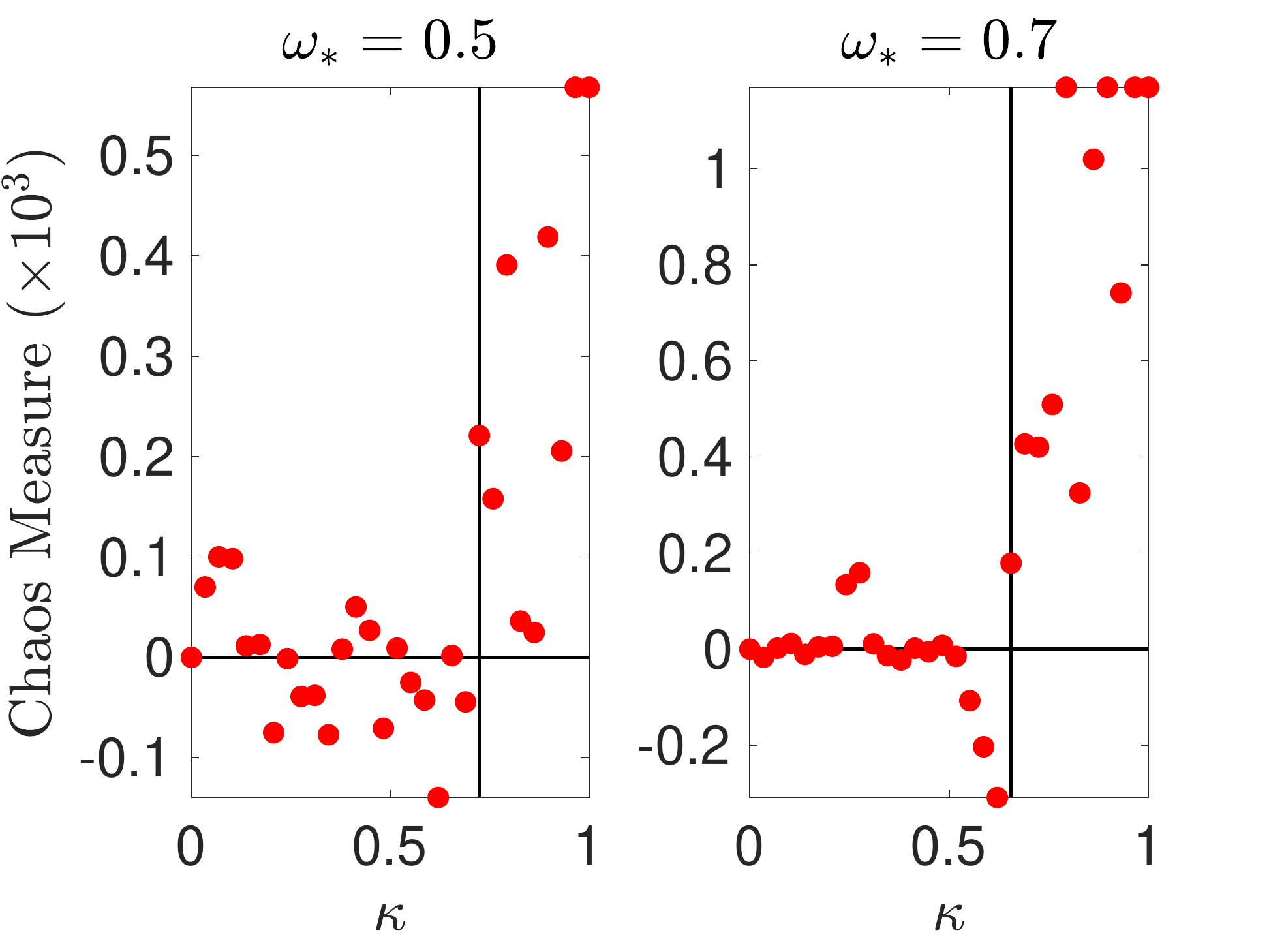}
\caption{\edit{Values of $\mathcal C$ as a function of $\kappa$ for the fixed frequencies $\omega_*=0.5$ and $\omega_*=0.7$.}}
\label{fig:chaos_ex}
\end{figure}

The chaotic interface (CI) is defined to be the set of pairs $(\omega_*,\kappa^*)$, where $\kappa^*(\omega_*)$ is the lower bound that satisfies $\mathcal C>0$ for all $\kappa>\kappa^*(\omega_*)$.  If no such value exists, we define $\kappa^*(\omega_*)=1$.  We say that the solution corresponding to $(\omega_*,\kappa)$ is chaotic if $\kappa>\kappa^*(\omega_*)$.  The value of $\kappa^*(\omega_*)$ is determined numerically from the simulations that produced the chaos measure in Fig. \ref{fig:chaos_measure}.  The methodology of computing $\kappa^*(\omega_*)$ is depicted in Fig. \ref{fig:chaos_ex}.  In this figure, $\mathcal C$ is plotted as a function of discrete amplitude factors ranging between $0\le \kappa\le 1$ for two fixed frequencies.  The horizontal and vertical lines described by $\mathcal C=0$ and $\kappa=\kappa^*(\omega_*)$ are depicted for reference.  It should be noted that small and positive values of $\mathcal C$ do not necessarily imply that the associated time series diverge.  This is to be expected because the flow variables need not oscillate sinusoidally and may favor more time near the peak of the wave.

\begin{figure}
\centering
\includegraphics[width=\linewidth]{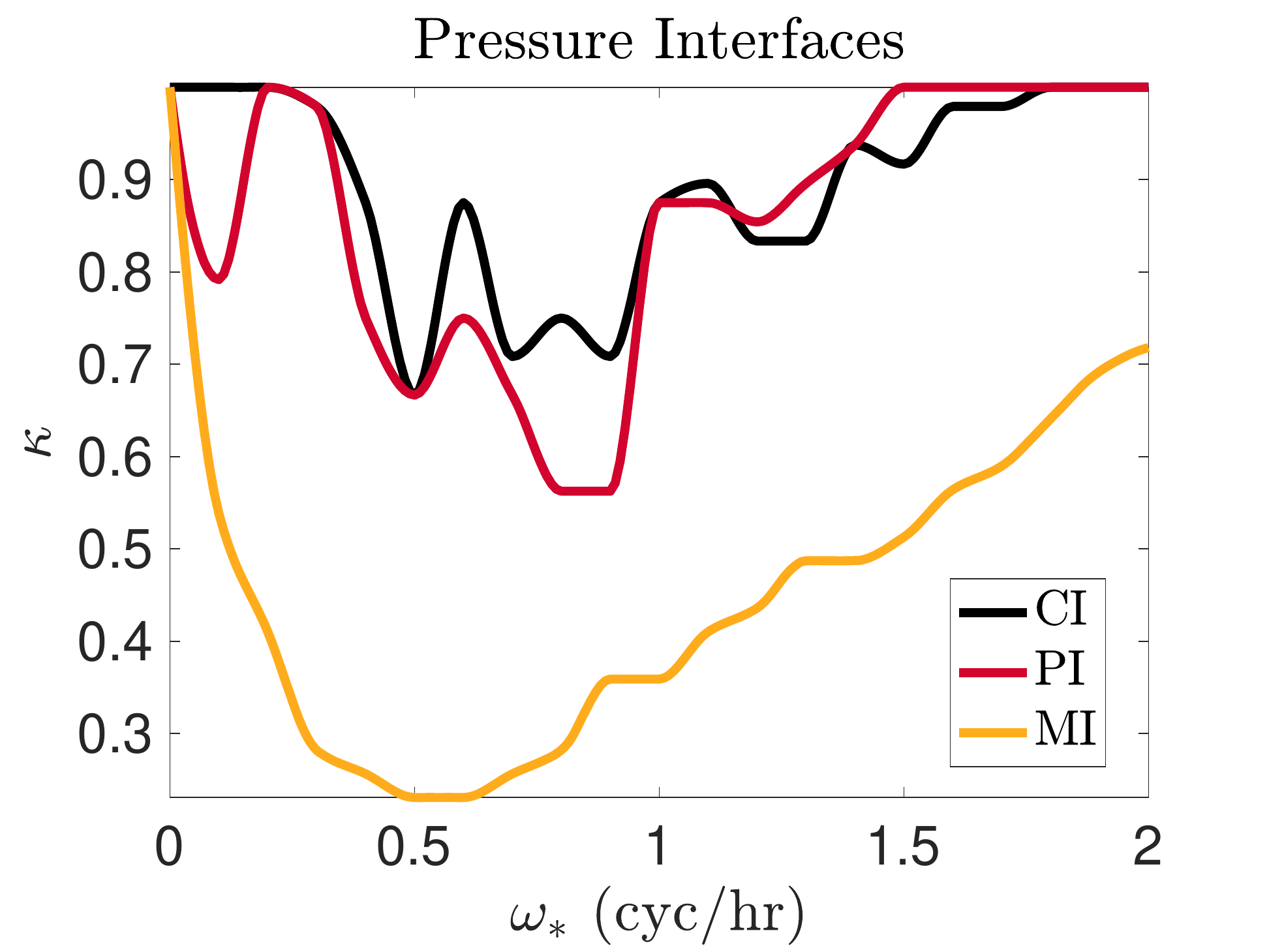}
\caption{\edit{Monotonic, periodic, and chaotic interfaces of the outlet pressure for the 50 (km) pipeline used in Fig. \ref{fig:powerspectrum} with $\alpha_1=0.2$.}}
\label{fig:Interfaces}
\end{figure}

We conclude our analysis with a comparison of the interfaces in the regions of boundary condition parameters $(\omega_*,\kappa)$ that do or do not exhibit monotonic, periodic, and chaotic responses to periodic boundary conditions for a single pipeline.  These interfaces, which we refer to as the monotonic interface (MI), periodic interface (PI), and chaotic interface (CI), are depicted in Fig. \ref{fig:Interfaces} for the single pipeline system used in Sections \ref{sec:MI} and \ref{sec:PI}.  The MI for pressure in Fig. \ref{fig:Interfaces} is different from the MI for pressure in Fig. \ref{fig:monotonic_interface} because of the different mean concentrations $\alpha_1$ used to compute the two MIs.  Moreover, the withdrawal rates of the three solutions that are used here for the computation of the current MI are $\bm w_2=40 \pi(D/2)^2$, $\bm w_2=75 \pi(D/2)^2$ and $\bm w_2=110 \pi(D/2)^2$, which are different from those used in Section \ref{sec:MI}.  There are a few key takeaways from Fig. \ref{fig:Interfaces}.  First, the region beneath the MI is a subset of the region beneath both the PI and CI and is significantly beneath these two interfaces for nonzero frequencies.  This suggests that all sinusoidal hydrogen blending conditions that show monotonicity properties will not result in irregular or chaotic responses, as expected for the strongly dissipative behavior of natural gas pipeline flows.  This result is informative for natural gas pipeline operations with significant hydrogen blending. Second, the PI and CI are qualitatively similar to one another.  Strictly speaking, the PI is beneath the CI for all $\omega_*\le 1.1$ and the CI is beneath the PI for all $\omega_*> 1.1$.  We note that the interfaces depend on the computational method and its parameter values, the discretization grid in the space of boundary condition parameters, and the threshold chosen to delimit the interfaces (i.e., $\mathcal P=0.3$ for the PI and $\mathcal C=0.0$ for the CI).  There may be other reasonable thresholds that yield somewhat different quantitative results, such as one where the PI is always beneath the CI or one where the PI is exactly equal to the CI at all of the grid points.  The qualitative similarity of the PI and CI suggests that these regions may indeed coincide, and only one of these measures may be sufficient to identify both non-periodic and chaotic behavior.
} 
 
\vspace{-2ex}

\section{Conclusions} \label{sec:conclusion}

We have developed a model for transporting heterogeneous mixtures of natural gas and hydrogen through pipeline networks.  The formulation may be applied to real pipeline systems with time-varying operations of compressor and regulator units, supply stations that inject gas into the network at defined pressure and hydrogen fraction, and flow stations that withdraw the mixture from the network.  The nonlinear partial differential equation formulation is discretized using a finite volume method to obtain a nonlinear input-to-state system, for which we prove monotone ordering properties for injections with constant hydrogen fraction, and prove that such monotonicity properties do not hold in general for injections with time-varying hydrogen concentrations.  This result builds on previous work that assumed globally homogeneous gas \cite{misra2020monotonicity}, by considering inhomogeneity in space (Proposition 1 on injections with constant hydrogen fraction) and time (Proposition 2 on injections with time-varying hydrogen fraction).  The interface in the boundary condition parameter region of concentration variation that partitions monotonic and non-monotonic responses was analyzed numerically and the results were illustrated on a test network.  Operations outside of the monotone response region may create surges with large pressure, energy, and concentration gradients, which do not occur in flows of a homogeneous gas.  The monotonic interface analysis indicates that sufficiently slow variation in concentration about a constant profile will likely maintain monotonicity of ordered solutions in overall system pressures, and prevent large, rapid pressure transients.  Such conditions are critical to maintain a physical flow regime with behavior that is intuitive for pipeline control room operators.   This suggests that hydrogen may be blended into a natural gas pipeline network as long as injection rates are changed only gradually.  The acceptable ramping rates depend significantly on the structure of the network, and would have to be determined through numerous simulations.

\edit{This study demonstrates that heterogeneous location and time dependent boundary conditions may result in non-periodic and chaotic flows when hydrogen is blended into a natural gas pipeline.  For a single pipeline, boundary conditions with monotonic, periodic, and chaotic responses were analyzed numerically and interfaces between regions where these properties do and do not hold were estimated.  The interface analysis of a single pipeline demonstrates that sinusoidal boundary conditions for which monotonicity properties are preserved lead to solutions that are not chaotic and will eventually approach a periodic response.  We have also demonstrated that the interfaces that delimit periodic and chaotic responses to periodic boundary conditions are qualitatively similar and may indeed coincide for appropriate thresholds that are used to define these interfaces.  This suggests that only one of these interfaces may be sufficient to identify both non-periodic and chaotic behavior.  Characterizing the monotonic, periodic, and chaotic interfaces with extensive simulations for specific network topologies may enable a gas pipeline system designer to determine limitations on operating their networks safely and predictably given blending of heterogeneous gases.}

\edit{The results of this study enable new capabilities to plan for and operate the transmission of gases through pipeline energy systems.  The generalizable modeling presented here can be used in addition to empirical studies to quantify the effect of hydrogen blending on gas pipelines under transient conditions \cite{styblo2022effects}.  The developed model can be applied in general natural gas operations if two or more production sites supply notably different levels of methane in their natural gas compositions.  More generally, the derivation of the gas mixture model may be extended from a mixture of two gases to any finite number of distinct gases.  We aim to present a self-contained analysis of a variety of different physical properties of gas mixture flows in pipeline networks.   Much of our analysis was performed on a finite-dimensional realization of the infinite-dimensional PDE system for a small network and a single pipeline.  Therefore, there are still a number of questions that remain open in terms of monotonicity, periodicity, and chaoticity.  Rigorous definitions and theorems regarding the limits of monotonicity, periodicity, and chaoticity of nonlinear PDEs on networks and their dependencies on graph topology have yet to be examined.  Such results may be used to establish a more complete analysis of appropriate blending conditions with which monotone-ordering and well-behaved flows are maintained. The onset of turbulence in complex multi-dimensional fluid mixing flows is extensively studied \cite{amon1996lagrangian,feigenbaum1980transition}, and here we show that complex, chaotic behavior can arise in a strongly dissipative, essentially one-dimensional distributed system. }

\clearpage


\section*{Acknowledgements}

The authors are grateful to Vitaliy Gyrya, Rodrigo Platte, Dieter Armbruster, and Yan Brodskyi for numerous helpful discussions, and to E. Olga Skowronek for drawing the networks in Figs. \ref{bc_conf} and \ref{net_conf}.  This study was supported by the U.S. Department of Energy's Advanced Grid Modeling (AGM) project ``Dynamical Modeling, Estimation, and Optimal Control of Electrical Grid-Natural Gas Transmission Systems'', as well as LANL Laboratory Directed R\&D project ``Efficient Multi-scale Modeling of Clean Hydrogen Blending in Large Natural Gas Pipelines to Reduce Carbon Emissions''. Research conducted at Los Alamos National Laboratory is done under the auspices of the National Nuclear Security Administration of the U.S. Department of Energy under Contract No. 89233218CNA000001.

\bibliography{sample}

\appendix 

\clearpage

\section{Model Comparison} \label{sec:comparison}

\begin{figure}[!th]
\centering
\includegraphics[width=\linewidth]{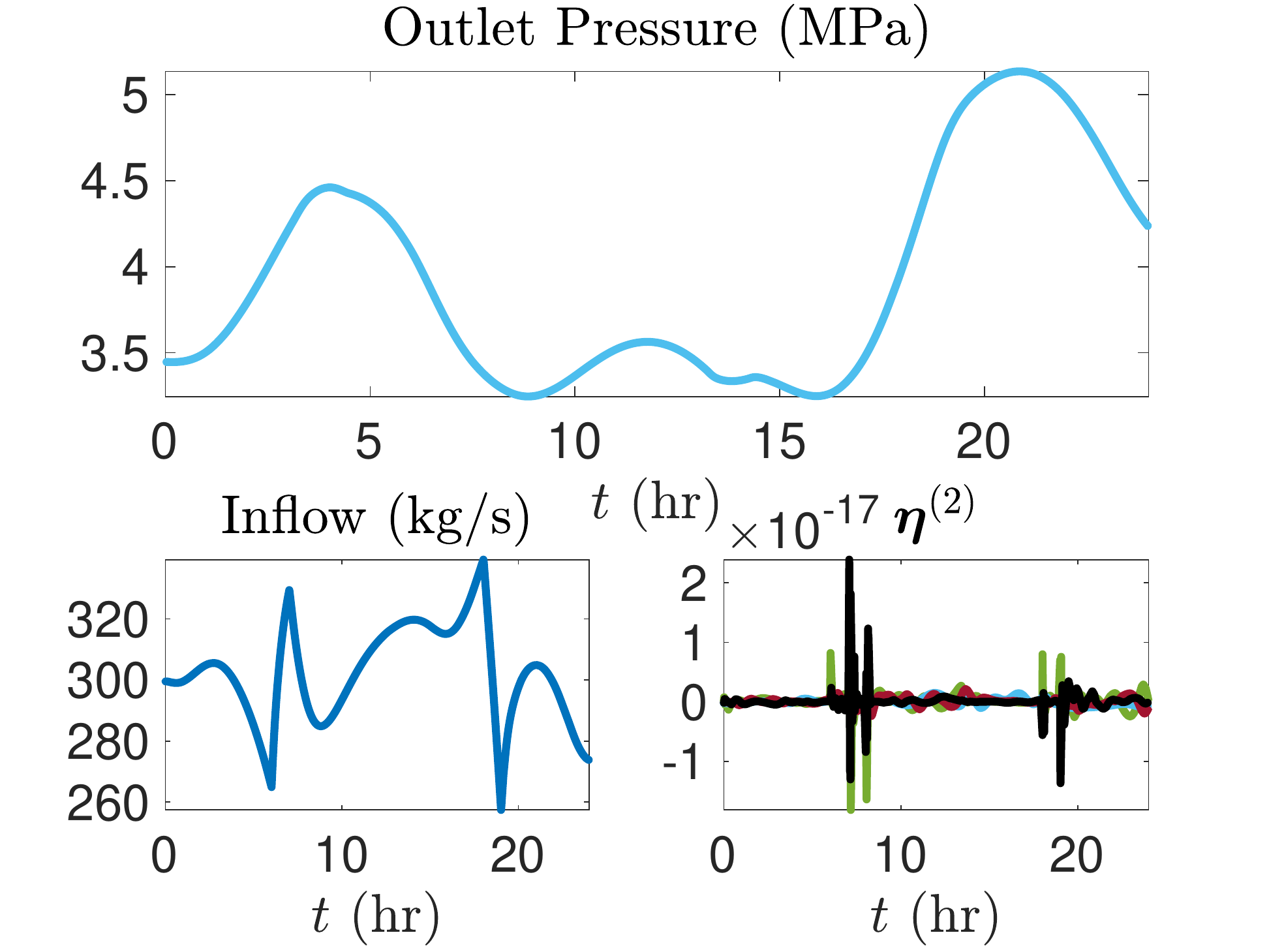}
\caption{Model comparison (see Fig. 11 in \cite{GYRYA201934}).}
\label{net_model_comp}
\end{figure}

The interested reader is referred to \cite{GYRYA201934} for a model comparison of homogeneous gas flow in the network shown in Fig. \ref{net_conf}.  Using the same initial and boundary conditions as used in the previous study, we recover the same solution, up to machine precision, with our mixed gas model.  Our solution is shown in Fig. \ref{net_model_comp}.

\section{Chebyshev Spectral Differentiation} \label{sec:chebyshev}

Consider a single pipeline of length $\ell$, diameter $D$, and friction factor $\lambda$ with axial variable $x\in [0,\ell]$.  Discretize the interval $[0,\ell]$ with the $(N+1)$ discretization points $x_i=\ell/2(1-\cos(i \pi/N))$ for $i=0,\dots,N$.  Define the sampled variables $\bm \rho^{(m)}_i(t)=\rho^{(m)}(t,x_i)$ and $\bm \varphi_i(t)=\varphi(t,x_i)$.  It follows from interpolating the values of $\bm \rho^{(m)}_i(t)$ at the points $x_i$ using Lagrange polynomials of order $N$ that (e.g., see \cite{ascher2011first})
\begin{equation}
    \partial_x \rho^{(m)}(t,x_i) \approx \bm D \bm \rho^{(m)}_i(t),
\end{equation}
where 
\begin{eqnarray}
\bm D_{ij} =
    \begin{cases}
    \sum \limits_{\substack{n=0 \\ n\neq j}}^n \frac{1}{x_j-x_n}, & i=j, \\
        \frac{1}{x_j-x_i}\prod \limits_{\substack{n=0 \\ n\neq i,j}}^n \frac{x_i-x_n}{x_j-x_n} & i\not= j.
    \end{cases}
\end{eqnarray}
The discretized PDEs in Eqs. \eqref{eq:pde1}-\eqref{eq:pde2} become
\begin{eqnarray}
 \dot{\bm \rho}^{(m)} +\bm D \left(\frac{ \bm \rho^{(m)}}{\bm \rho^{(1)}+\bm \rho^{(2)}} \odot \bm \varphi \right) &=&0, \label{eq:cheb1} \\
\bm D \left(\sigma^2_1 \bm \rho^{(1)}+\sigma^2_2 \bm \rho^{(2)} \right) &=& -\frac{\lambda}{2D}\frac{ \bm \varphi \odot |  \bm \varphi|}{\bm \rho^{(1)}+\bm \rho^{(2)}}. \label{eq:cheb2} \quad
\end{eqnarray}
The boundary conditions are incorporated into the discretized equations by replacing $\bm \rho^{(m)}_0(t)=\bm s_0^{(m)}(t)$ and $\bm \varphi_N(t)=\bm w_N(t)/(0.25\pi D^2)$.

\edit{
\section{Nomenclature} \label{sec:nomenclature}

\begingroup
\addtolength{\jot}{-1pt}
 \noindent {\bf Network Graph}
    \begin{align*}
  \mathcal E &     && \text{Directed Edge Set} \\
   \mathcal V &     && \text{Node Set} \\
   \mathcal V_s \subset \mathcal V &     && \text{Slack Nodes} \\
    \mathcal V_d \subset \mathcal V &     && \text{Non-slack Nodes} \\
    \mathcal V_w \subset \mathcal V_d &     && \text{Non-slack Withdrawal Nodes} \\
    \mathcal V_q \subset \mathcal V_d &     && \text{Non-slack Injection Nodes} \\
      E &     && \text{Cardinality of $\mathcal E$} \\
      V &     && \text{Cardinality of $\mathcal V$} \\
      V_s  &     && \text{Cardinality of $\mathcal V_s$} \\
      V_d  &     && \text{Cardinality of $\mathcal V_d$} \\
      V_w  &     && \text{Cardinality of $\mathcal V_w$} \\
    V_q &     && \text{Cardinality of $\mathcal V_q$} \\
   k &     && \text{Edge Index} \\
   i \; \& \; j &     && \text{Node Indices} \\
    k:i\mapsto j &     && \text{Edge $k$ Directs from Node $i$ to Node $j$} \\
    _{\mapsto}j &     && \text{Set of Edges Directed to Node $j$} \\
     j_{\mapsto} &     && \text{Set of Edges Directed from Node $j$} 
     \end{align*}
\noindent {\bf Edge Variables} 
    \begin{align*}
   \rho &     && \text{Total Density (kg m$^{-3}$)} \\
    p &     && \text{Total Pressure (MPa)} \\
    \varphi &     && \text{Total Mass Flux (kg m$^{-2}$ s$^{-1}$) } \\
     \rho^{(m)} &     && \text{Partial Density (kg m$^{-3}$)} \\
     p^{(m)} &     && \text{Partial Pressure (MPa)} \\
     \varphi^{(m)} &     && \text{Partial Mass Flux (kg m$^{-2}$ s$^{-1}$) } \\
     \eta^{(m)} &     && \text{Mass Fraction} \\
     \nu^{(m)} &     && \text{Volumetric Fraction} \\
       \underline \psi &     && \text{Evaluation of $\psi$ at Edge Inlet} \\
  \overline \psi &     && \text{Evaluation of $\psi$ at Edge Outlet} 
            \end{align*}
\noindent {\bf Node Variables} 
  \begin{align*}
      \bm \rho &     && \text{Total Density (kg m$^{-3}$)} \\
    \bm p &     && \text{Total Pressure (MPa)} \\
            \bm E &     && \text{Energy (GJ)} \\
    \bm \rho^{(m)} &     && \text{Partial Density (kg m$^{-3}$)} \\
     \bm \eta^{(m)} &     && \text{Mass Fraction} \\
     \bm \nu^{(m)} &     && \text{Volumetric Fraction} 
     \end{align*}
\noindent {\bf Boundary Condition Variables} 
\begin{align*}
 \bm s^{(m)} &     && \text{Slack Node Partial Density (kg m$^{-3}$)} \\
  \bm p_{s} &     && \text{Slack Node Pressure (MPa)} \\
  \bm \alpha &     && \text{Slack Node Mass Fraction} \\
   \bm \beta &     && \text{Non-slack Injection Node Mass Fraction} \\
     \bm a &     && \text{Slack Node Wave Speed} \\
   \bm b &     && \text{Non-slack Injection Node Wave Speed} \\
   \bm w &     && \text{Non-slack Withdrawal Node Mass Outflow (kg s$^{-1}$)} \\
    \bm q &     && \text{Non-slack Injection Node Mass Inflow (kg s$^{-1}$)} 
    \end{align*}
\noindent {\bf Control Variables} 
\begin{align*}
  \underline \mu &     && \text{Compressor Ratio} \\
  \overline \mu &     && \text{Regulator Ratio} 
      \end{align*}
\noindent {\bf Network Graph Matrices}
\begin{align*}
  M &     && \text{Weighted Incidence Matrix} \\
   M_s &     && \text{Weighted Slack Node Incidence Submatrix} \\
   M_d &     && \text{Weighted Non-slack Node Incidence Submatrix} \\
    Q &     && \text{Incidence Matrix} \\
    Q_s &     && \text{Slack Node Incidence Submatrix} \\
   Q_d &     && \text{Non-slack Node Incidence Submatrix} \\   
    \underline A &     && \text{Componentwise Negative Parts of $A$} \\
    \overline A &     && \text{Componentwise Positive Parts of $A$} 
      \end{align*}
 \noindent {\bf Parameters}     
      \begin{align*}
     \sigma_1 &     && \text{Natural Gas Wave Speed} \\
       \sigma_2 &     && \text{Hydrogen Wave Speed} \\
       \sigma &     && \text{Local Mixture Wave Speed} \\
       \ell &     && \text{Pipe Length (km)} \\
       L &     && \text{Diag($\ell_k$)} \\
        D &     && \text{Pipe Diameter (m)} \\
        \chi &     && \text{Pipe Cross-sectional Area (m$^2$)} \\
        X &     && \text{Diag($\chi_k$)} \\
       \lambda &     && \text{Darcy-Weisbach Friction Factor}  \\
       \Lambda &     && \text{Diag($\sqrt{2D_k/\lambda_k\ell_k}$)}  
      \end{align*}
\endgroup
}

\clearpage

\end{document}